\documentclass[pdflatex,sn-mathphys-num]{sn-jnl}

\usepackage{amsmath,amssymb,amsfonts}
\usepackage{amsthm}
\usepackage{booktabs}
\usepackage{graphicx}
\usepackage{mathtools}

\graphicspath{{./}}

\numberwithin{equation}{section}

\theoremstyle{thmstyleone}
\newtheorem{theorem}{Theorem}
\newtheorem{corollary}[theorem]{Corollary}
\newtheorem{proposition}[theorem]{Proposition}
\newtheorem{lemma}[theorem]{Lemma}

\theoremstyle{thmstyletwo}
\newtheorem{example}{Example}
\newtheorem{remark}{Remark}

\begin{document}

\title[Bilevel Optimization for CLG Optical Flow]{A Theoretical Analysis of a Bilevel Optimization Formulation for Combined Local-Global Optical Flow Models}

\author*[1]{\fnm{El Mostafa} \sur{Kalmoun}}\email{E.Kalmoun@aui.ma}

\affil*[1]{\orgdiv{School of Science and Engineering}, \orgname{Al Akhawayn University in Ifrane},
\orgaddress{\city{Ifrane}, \postcode{53000}, \country{Morocco}}}

\abstract{This paper analyzes a bilevel formulation for learning the regularization parameter
of combined local-global (CLG) optical flow models. The lower-level problem computes the flow from a
Lucas--Kanade data term and a quadratic or total variation (TV) regularizer, and the upper-level problem
selects the parameter that minimizes the squared distance to a ground-truth flow. The setting is
finite-dimensional, and a full-rank assumption makes the lower-level solution unique. For the quadratic
model we derive sufficient conditions for the existence of an optimal parameter, for every optimal
parameter to be positive, and for strict local optimality. For TV regularization an optimal parameter
always exists, because the lower-level solution vanishes for large parameters, and we give a sufficient
condition for positivity that accounts for the coupling of the two flow components through the local
structure tensors. Synthetic examples illustrate the conditions, including cases in which they fail.}

\keywords{Bilevel Optimization, Optical Flow, Regularization Parameter Learning, Total Variation, Local-Global Models}

\maketitle

\section{Introduction}\label{sec:intro}

Variational methods estimate the optical flow between two frames by minimizing an energy that combines
a data term, derived from brightness constancy, with a regularizer weighted by a parameter $\alpha>0$. The
quality of the estimate depends strongly on $\alpha$. Too small a value leaves the flow noisy, and too
large a value oversmooths it. This paper studies the choice of $\alpha$ by bilevel learning from
ground-truth flow and asks, in a fully discrete setting, when an optimal parameter exists and when it is
positive.

The linearized brightness constancy
constraint gives one equation per pixel for the two flow components, so the flow is not determined
locally (the aperture problem). Lucas and Kanade~\cite{lucas1981iterative} assume that the flow is constant
in a window around each pixel and solve the resulting local least-squares problems. Horn and
Schunck~\cite{horn1981determining} instead add a global quadratic penalty on the flow gradient, which
yields dense flow fields but smooths motion boundaries. Total variation (TV)
regularization, introduced for image denoising~\cite{rudin1992nonlinear}, preserves such boundaries.
Shulman and Herv\'e~\cite{shulman1989regularization} proposed an early non-quadratic regularizer for
discontinuous flow fields, and TV is widely used for optical
flow~\cite{rouchouze1994motion,deriche1995optical,zach2007aduality,Wedel2008,chambolle2011afirst,rantfl2014nonlocal,ali2016illumination}.
We use its isotropic form, applied to each flow component. TV is convex but not differentiable, and it is
usually minimized by primal--dual methods~\cite{chambolle2004algorithm,zach2007aduality,chambolle2011afirst} or
replaced by a smooth approximation such as the Huber function~\cite{kalmoun2018investigation}. The combined
local-global (CLG) model of Bruhn et al.~\cite{bruhn2005lucas} joins the Lucas--Kanade data term with a global
regularizer. It inherits the robustness of the local method and the density of the global one, and it
can be minimized efficiently, for instance by multigrid methods~\cite{kalmoun2013trust}. Spatially varying
parameters $\alpha(p)$ adapted to image features have also been
proposed~\cite{wedel2009structure,Belhachmi2011,Belhachmi2016}.

Given two frames and a ground-truth flow $w_g\in\mathbb{R}^{2N}$ on a grid $\Omega$ of $N$ pixels, we consider
\begin{equation}\label{eq:intro_problem}
	\min_{\alpha\ge0}\ \|w_\alpha-w_g\|_2^2
	\quad\text{subject to}\quad
	w_\alpha=\operatorname*{arg\,min}_{w\in\mathbb{R}^{2N}}\ \frac12\sum_{p\in\Omega}\|A_pw(p)-b_p\|_2^2+\alpha R(w),
\end{equation}
where $A_pw(p)=b_p$ is the Lucas--Kanade system of the window around the pixel $p$, and $R$ is either a
quadratic penalty on the discrete flow gradient or the TV seminorm. Section~\ref{sec:setting} defines all
objects precisely.

Problem~\eqref{eq:intro_problem} rests on four modeling choices. First, the parameter is learned from
ground truth. Classical parameter-choice rules, such as the discrepancy
principle~\cite{engl1996regularization}, the L-curve~\cite{hansen1992lcurve}, generalized
cross-validation~\cite{golub1979gcv} and Stein's unbiased risk estimate~\cite{stein1981sure}, need no ground
truth. They rely instead on knowledge of the noise or on heuristic criteria, and they do not target the
flow error itself. Bilevel learning uses training data with known ground truth and minimizes the
reconstruction error directly. Second, the upper-level loss is the squared
Euclidean distance, which equals the sum over all pixels of the squared endpoint errors, a standard
accuracy measure for optical flow~\cite{baker2011middlebury}. It is smooth in $w$, whereas the endpoint error
itself and the angular error~\cite{barron1994performance} are not differentiable where the error vanishes.
Third, a single scalar parameter is the natural first case. It is the setting in which existence and
positivity can be characterized by explicit conditions, and spatially varying parameters, as in the
bilevel framework of~\cite{pagliari2022tvtype}, can build on it. Fourth, the CLG data term makes the
lower-level problem well posed. When every window contains two pixels with linearly independent
gradients, the data term is strongly convex, and the lower-level solution is unique for every
$\alpha\ge0$, including $\alpha=0$ (Lemma~\ref{lem:wellposed}). The Horn--Schunck data term alone is never
strongly convex, because its matrix $J^\top J$ in~\eqref{OFCE} has a null space of dimension at least~$N$.

Bilevel learning of regularization parameters in imaging was established by Kunisch and
Pock~\cite{kunisch2013bilevel} for variational denoising and by De los Reyes and
Sch{\"o}nlieb~\cite{delosreyes2013} for learning noise models in TV denoising. Its theory
includes existence results for a broad class of regularizers~\cite{delosreyes2016structure}, refined
stationarity systems for TV~\cite{delosreyes2022optimality}, stationarity and second-order conditions through
a reformulation as a mathematical program with complementarity constraints~\cite{delosreyes2023mpcc}, and
spatially dependent parameters~\cite{pagliari2022tvtype}. For finite-dimensional linear inverse problems
with a continuously differentiable regularizer, Ehrhardt, Gazzola and Scott~\cite{ehrhardt2024optimal} give
a sufficient condition for a positive optimal parameter. Our positivity result for the quadratic model
(Theorem~\ref{th:clg-lossE2}) is a special case of their condition. Applications include blind
deconvolution~\cite{hintermuller2015bilevel} and the learning of fractional differential
operators~\cite{bartels2020fractional}, and Crockett and Fessler~\cite{crockett2022bilevel} survey the field.
Computational work addresses inexact lower-level solves through derivative-free
methods~\cite{ehrhardt2021dfbilevel}, inexact hypergradients~\cite{ehrhardt2024hypergradients} and single-loop
schemes~\cite{suonpera2025singleloop}. Learning through truncated or approximate
inference~\cite{domke2012generic,domke2013approxmarg} and learned variational inference layers~\cite{ranftl2014}
aim at empirical performance rather than at the analysis of a regularization weight. For optical flow,
Jiang et al.~\cite{jiang2020joint} use a bilevel scheme that couples flow and camera motion in unsupervised
deep learning. Closest to the present setting, El Hakoume et al.~\cite{elhakoume2023lame} learn the two
Lam\'e coefficients of a fluid optical flow registration model by bilevel supervised learning, within a
multi-frame super-resolution method. Positivity is imposed through a fixed box, uniqueness of the
lower-level solution is assumed, and the bilevel problem ``is not theoretically studied'' because their
model lacks the required regularity. Related work by Laghrib and coauthors learns spatially varying weights
for nonlocal deblurring~\cite{elmalki2024nonlocal}. To our knowledge, no previous work establishes existence,
positivity and second-order conditions for learning the regularization parameter of a variational
optical flow model from ground-truth flow. Here all
objects are finite-dimensional, the proofs are direct, and the conditions are explicit and computable from
the data for both regularizers.

The main contributions are the following.
\begin{itemize}
	\item A finite-dimensional bilevel formulation for CLG optical flow in which a scalar regularization
	parameter is learned from ground-truth flow.
	\item Existence of an optimal parameter. For any convex regularizer that vanishes only at zero, a
	sufficient condition is that the Lucas--Kanade estimate be closer to the ground truth than the zero
	flow (Proposition~\ref{prop:exist}). For TV with zero boundary values this condition is not needed,
	because the lower-level solution vanishes beyond a finite value of $\alpha$ (Corollary~\ref{thm:tv_exist}).
	\item Sufficient conditions under which $\alpha=0$ is not optimal, so that every optimal parameter is
	positive, for the quadratic model
	(Theorem~\ref{th:clg-lossE2}, an instance of the condition of~\cite{ehrhardt2024optimal}) and for TV
	(Theorem~\ref{thm:tv_nonzero}), and second-order sufficient conditions for strict local optimality in the
	quadratic model (Theorem~\ref{th:sosc}). The TV condition is not covered by~\cite{ehrhardt2024optimal},
	because TV is not differentiable, and it differs structurally from the analyses of scalar data
	terms~\cite{delosreyes2016structure,delosreyes2022optimality}.
	The Lucas--Kanade data term couples the two flow components through $2\times2$ blocks that vary from
	pixel to pixel. This coupling requires a four-way partition of the pixels, a term
	$\mathcal{R}_{\mathrm{mix}}$ for the pixels at which only one flow component has a nonzero gradient, and a
	term $\mathcal{R}_{\mathrm{var}}$ for the spatial variation of the blocks. These terms vanish for spatially
	uniform scalar data terms.
	\item Synthetic examples that illustrate the conditions, including cases in which the sufficient
	conditions fail.
\end{itemize}

Section~\ref{sec:setting} introduces the discrete setting and states the bilevel problem.
Section~\ref{sec:quad} analyzes quadratic regularization and Section~\ref{sec:TV} total variation
regularization. Section~\ref{sec:examples} presents the examples, and Section~\ref{sec:conclusion} discusses
the limitations of the analysis and concludes.

\section{Discrete Setting and Problem Formulation}\label{sec:setting}

All objects in this paper are finite-dimensional. The image domain $\Omega$ is a rectangular grid of
$N=n_yn_x$ pixels with unit grid spacing, ordered row by row. Scalar fields on $\Omega$, such as images,
are vectors in $\mathbb{R}^N$. A flow field is $w=(u;v)\in\mathbb{R}^{2N}$ with $u,v\in\mathbb{R}^N$, and its
value at a pixel $p$ is $w(p)=(u(p),v(p))^\top\in\mathbb{R}^2$.

For two frames $I_1,I_2\in\mathbb{R}^N$, let $I_x,I_y\in\mathbb{R}^N$ be the spatial derivatives and $I_t:=I_2-I_1$
the temporal derivative. The computation of $I_x$ and $I_y$ is described at the beginning of
Section~\ref{sec:examples}. With
$\nabla I(p):=(I_x(p),I_y(p))^\top$ and
$J:=\bigl[\operatorname{diag}(I_x)\ \ \operatorname{diag}(I_y)\bigr]\in\mathbb{R}^{N\times2N}$, the linearized
brightness constancy constraint reads
\begin{equation}\label{OFCE}
	Jw+I_t=0\in\mathbb{R}^N,\qquad\text{that is,}\qquad\nabla I(p)^\top w(p)+I_t(p)=0\quad\text{for every }p\in\Omega .
\end{equation}

The Lucas--Kanade data term assigns to each pixel $p$ a window $\Omega_p\subset\Omega$ that contains $p$,
and $n_p:=|\Omega_p|$. Imposing~\eqref{OFCE} at every $q\in\Omega_p$ with the common
value $w(p)$ gives the local system $A_pw(p)=b_p$ with
\begin{equation}\label{LK-A-b}
	A_p:=\bigl[\nabla I(q)^\top\bigr]_{q\in\Omega_p}\in\mathbb{R}^{n_p\times2},\qquad
	b_p:=-\bigl[I_t(q)\bigr]_{q\in\Omega_p}\in\mathbb{R}^{n_p}.
\end{equation}
The data term is
\begin{equation}\label{eq:data_term}
	D(w):=\frac12\sum_{p\in\Omega}\|A_pw(p)-b_p\|_2^2=\frac12\,w^\top\mathcal{M}_{\mathrm{LK}}w-f_{\mathrm{LK}}^\top w+c,
\end{equation}
where $f_{\mathrm{LK}}\in\mathbb{R}^{2N}$ collects the vectors $A_p^\top b_p$ and $c:=\frac12\sum_{p\in\Omega}\|b_p\|_2^2$. Let $\Pi$ be the permutation matrix that orders the unknowns pixel by pixel:
\[\Pi w=(u(p_1),v(p_1),\dots,u(p_N),v(p_N))^\top.\] 
Stacking all local systems gives the global system
\begin{equation}\label{LK}
	Aw=b,
\end{equation}
whose matrix $A$ has $n:=\sum_pn_p$ rows and $2N$ columns, and $A\Pi^\top=\operatorname{diag}(A_{p_1},\dots,A_{p_N})$ is block
diagonal with rectangular $n_p\times2$ blocks. Hence
\[
	\mathcal{M}_{\mathrm{LK}}=A^\top A=\Pi^\top\operatorname{diag}\bigl(A_{p_1}^\top A_{p_1},\dots,A_{p_N}^\top A_{p_N}\bigr)\,\Pi .
\]
In the ordering $w=(u;v)$, $\mathcal{M}_{\mathrm{LK}}$ consists of four diagonal $N\times N$ blocks formed
from the entries of the matrices $A_p^\top A_p$.

Throughout the paper we assume
\begin{center}
	\textbf{(A1)}\quad $\operatorname{rank}A_p=2$ for every $p\in\Omega$,
\end{center}
that is, every window contains two pixels with linearly independent image gradients. Then
$A_p^\top A_p\succ0$ for every $p$, so $\mathcal{M}_{\mathrm{LK}}\succ0$, and $D$ has the unique minimizer
$w_0=\mathcal{M}_{\mathrm{LK}}^{-1}f_{\mathrm{LK}}$, the Lucas--Kanade solution
\begin{equation}\label{LSS}
	\begin{gathered}
		w_0(p)=(A_p^\top A_p)^{-1}A_p^\top b_p,\\
		A_p^\top A_p=\sum_{q\in\Omega_p}\nabla I(q)\nabla I(q)^\top,\qquad
		A_p^\top b_p=-\sum_{q\in\Omega_p}I_t(q)\,\nabla I(q).
	\end{gathered}
\end{equation}

The two regularizers use two difference operators, and their boundary conditions are part of the
model. For the quadratic regularizer, $K_D$ takes forward differences between horizontally and vertically
adjacent pixels, with zero ghost values on all four sides of the grid, so $K_D\in\mathbb{R}^{m\times N}$ with
$m=n_y(n_x+1)+n_x(n_y+1)$. The matrix
$K_D^\top K_D$ is the five-point Laplacian with homogeneous Dirichlet conditions. On flows we set $\mathbf{K}_D:=\operatorname{diag}(K_D,K_D)$ and
\[
	L:=\mathbf{K}_D^\top\mathbf{K}_D=\operatorname{diag}\bigl(K_D^\top K_D,\,K_D^\top K_D\bigr)\in\mathbb{R}^{2N\times2N}.
\]
Since $K_D$ is injective, $L\succ0$. The matrix $\mathcal{L}_{\mathrm{glob}}:=-L$ is the discrete vector
Laplacian.

For the TV regularizer, the gradient of a scalar field $s$ at a pixel $p$ is
\[
	\nabla s(p):=(Ks)(p):=\bigl(\bar s(p+e_x)-s(p),\ \bar s(p+e_y)-s(p)\bigr)^\top\in\mathbb{R}^2,
\]
where $e_x$ and $e_y$ are the shifts to the next column and the next row, and $\bar s$ equals $s$ on
$\Omega$ and $0$ outside $\Omega$. Thus the ghost values beyond the last column and the last row are
zero. The map $K:\mathbb{R}^N\to\mathbb{R}^{2N}$ is injective. On flows,
$\mathbf{K}:=\operatorname{diag}(K,K):\mathbb{R}^{2N}\to\mathbb{R}^{4N}$, and a dual field is
$z=(z_u,z_v)\in\mathbb{R}^{4N}$ with $z_u(p),z_v(p)\in\mathbb{R}^2$. The discrete divergence is
$\operatorname{div}:=-K^\top$ for scalar fields and $\mathbf{div}:=-\mathbf{K}^\top:\mathbb{R}^{4N}\to\mathbb{R}^{2N}$
for flows, so that
\begin{equation}\label{eq:adjoint}
	\sum_{p\in\Omega}\langle\nabla s(p),z(p)\rangle=-\sum_{p\in\Omega}s(p)\,(\operatorname{div}z)(p)
\end{equation}
holds exactly, with no boundary terms. Forward differences of other fields on $\Omega$, such as
residuals, are also taken with $K$. The two operators differ at the border, since $K$ has no
differences across the left and the top sides, and $K^\top K\neq K_D^\top K_D$. The analysis of the
quadratic model uses only $L\succ0$. For TV, existence uses only the injectivity of $K$, and the
positivity analysis also uses the forward-difference form of $K$ (Lemma~\ref{lem:prod}).

With these operators, the quadratic (CLG) energy and the TV energy are
\begin{align}
	E(w,\alpha)&:=D(w)+\frac{\alpha}{2}\|\mathbf{K}_Dw\|_2^2=D(w)+\frac{\alpha}{2}\,w^\top Lw,\label{eq:E_quad}\\
	E(w,\alpha)&:=D(w)+\alpha\,\mathcal{R}_{TV}(w),\qquad
	\mathcal{R}_{TV}(w):=\sum_{p\in\Omega}\bigl(\|\nabla u(p)\|_2+\|\nabla v(p)\|_2\bigr).\label{eq:E_TV}
\end{align}
The Horn--Schunck energy~\cite{horn1981determining}
\begin{equation}\label{HS}
	\tfrac12\|Jw+I_t\|_2^2+\tfrac{\alpha}{2}\|\mathbf{K}_Dw\|_2^2
\end{equation}
is the special case $\Omega_p=\{p\}$ of~\eqref{eq:E_quad}, for which (A1) fails.

Given a ground-truth flow $w_g\in\mathbb{R}^{2N}$, we consider
\begin{equation}\label{eq:bilevel_opt_flow}
	\min_{\alpha\ge0}\ Q(\alpha):=\|w_\alpha-w_g\|_2^2
	\quad\text{subject to}\quad
	w_\alpha=\operatorname*{arg\,min}_{w\in\mathbb{R}^{2N}}E(w,\alpha),
\end{equation}
with $E$ given by~\eqref{eq:E_quad} in Section~\ref{sec:quad} and by~\eqref{eq:E_TV} in
Section~\ref{sec:TV}. 

\begin{lemma}\label{lem:wellposed}
	Assume (A1), and let $R:\mathbb{R}^{2N}\to[0,\infty)$ be convex. For every $\alpha\ge0$, the function
	$D+\alpha R$ is strongly convex with modulus $\mu:=\lambda_{\min}(\mathcal{M}_{\mathrm{LK}})>0$ and has
	exactly one minimizer. In particular, $w_\alpha$ in~\eqref{eq:bilevel_opt_flow} is uniquely defined for
	both energies and every $\alpha\ge0$, including $\alpha=0$, and $Q$ is a single-valued function on
	$[0,\infty)$.
\end{lemma}

\begin{proof}
	Since $\mathcal{M}_{\mathrm{LK}}-\mu I\succeq0$, the function $D-\frac{\mu}{2}\|\cdot\|_2^2$ is convex, so
	$D+\alpha R$ is strongly convex. A finite strongly convex function on $\mathbb{R}^{2N}$ has exactly one
	minimizer.
\end{proof}

Without (A1), the lower-level problem may have many solutions, for instance at $\alpha=0$ because of
the aperture problem. Problem~\eqref{eq:bilevel_opt_flow} would then require an optimistic or a
pessimistic formulation, which we do not consider.

\begin{remark}\label{rem:weighted}
	When the ground truth is only partially available or carries confidence weights, the loss can be
	replaced by $Q_W(\alpha):=\|W(w_\alpha-w_g)\|_2^2$ with $W:=\operatorname{diag}(\omega,\omega)\in\mathbb{R}^{2N\times2N}$
	and weights $\omega\in\mathbb{R}^N$, $\omega\ge0$. The lower-level solutions do not depend on $W$. The
	corresponding changes to the conditions of Sections~\ref{sec:quad} and~\ref{sec:TV} are stated in
	Remark~\ref{rem:weighted_conditions}.
\end{remark}

\section{Quadratic Regularization}\label{sec:quad}

In this section the lower-level energy is the quadratic CLG energy~\eqref{eq:E_quad},
\begin{equation}\label{eq:CLG_expanded}
	E(w,\alpha)=\frac12\sum_{p\in\Omega}\|A_pw(p)-b_p\|_2^2+\frac{\alpha}{2}\,w^\top Lw .
\end{equation}
At the minimizer the gradient of $E(\cdot,\alpha)$ vanishes. For each pixel this gives
\begin{equation}\label{eq:EL_CLG_p}
	A_p^\top\bigl(A_pw(p)-b_p\bigr)+\alpha\,(Lw)(p)=0,
\end{equation}
where $(Lw)(p)\in\mathbb{R}^2$ collects the two entries of $Lw$ at $p$. Globally,
\begin{equation}\label{eq:EL_CLG_global}
	(\mathcal{M}_{\mathrm{LK}}+\alpha L)\,w=f_{\mathrm{LK}} .
\end{equation}
Under (A1), $\mathcal{M}_{\mathrm{LK}}\succ0$, and $L\succ0$. The matrix $\mathcal{M}_{\mathrm{LK}}+\alpha L$ is
therefore positive definite for every $\alpha\ge0$. It remains positive definite for small negative
$\alpha$: since $w^\top(\mathcal{M}_{\mathrm{LK}}+\alpha L)w\ge\bigl(\lambda_{\min}(\mathcal{M}_{\mathrm{LK}})
-|\alpha|\,\lambda_{\max}(L)\bigr)\|w\|_2^2$, this holds for every $\alpha>-\alpha_0$ with
$\alpha_0:=\lambda_{\min}(\mathcal{M}_{\mathrm{LK}})/\lambda_{\max}(L)>0$. Hence
\begin{equation}\label{eq:w_alpha_global}
	w_\alpha=X(\alpha)f_{\mathrm{LK}},\qquad X(\alpha):=(\mathcal{M}_{\mathrm{LK}}+\alpha L)^{-1},
\end{equation}
and $\alpha\mapsto w_\alpha$ and $Q$ are infinitely differentiable on the open interval $(-\alpha_0,\infty)$. In particular, $Q'(0)$ is an ordinary two-sided derivative, although only
$\alpha\ge0$ is admissible in the bilevel problem. From $X'(\alpha)=-X(\alpha)LX(\alpha)$ and $X(\alpha)Lw_\alpha=-w'_\alpha$,
\begin{equation}\label{eq:w_derivative}
	w'_\alpha=-X(\alpha)Lw_\alpha,\qquad w''_\alpha=-2X(\alpha)Lw'_\alpha,
\end{equation}
and therefore
\begin{equation}\label{eq:Q_gradient}
	Q'(\alpha)=2\langle w'_\alpha,\,w_\alpha-w_g\rangle,\qquad
	Q''(\alpha)=2\|w'_\alpha\|_2^2-4\langle X(\alpha)Lw'_\alpha,\,w_\alpha-w_g\rangle .
\end{equation}
Substituting~\eqref{eq:w_alpha_global} reduces~\eqref{eq:bilevel_opt_flow} to the single-level problem
\begin{equation}\label{eq:SL_CLG}
	\min_{\alpha\ge0}\ Q(\alpha)=\|X(\alpha)f_{\mathrm{LK}}-w_g\|_2^2 .
\end{equation}
If every minimizer satisfies $\alpha^*>0$, then some amount of regularization strictly improves on
the unregularized Lucas--Kanade estimate. 

\begin{remark}\label{rem:numerical_impl}
	In practice $X(\alpha)$ is never formed. The symmetric positive definite system~\eqref{eq:EL_CLG_global}
	is solved by preconditioned conjugate gradients or multigrid, and $w'_\alpha$ requires one more solve
	with the same matrix. Since $\alpha$ is a scalar, $Q$ can be minimized by a bounded one-dimensional
	search, which is what the examples of Section~\ref{sec:examples} use, or by a gradient-based method
	with the analytical gradient~\eqref{eq:Q_gradient}.
\end{remark}

The proof of the next existence result uses the minimizing-sequence argument
of~\cite[Proposition~3.1]{kunisch2013bilevel} in finite dimension.

\begin{proposition}\label{prop:exist}
	Assume (A1). Let $R:\mathbb{R}^{2N}\to[0,\infty)$ be convex with $R(w)=0$ if and only if $w=0$, and let
	$w_\alpha$ be the minimizer of $D+\alpha R$. Then:
	\begin{enumerate}
		\item[(i)] $\alpha\mapsto w_\alpha$ is continuous on $[0,\infty)$,
		$\alpha R(w_\alpha)\le\frac12w_0^\top\mathcal{M}_{\mathrm{LK}}w_0$ for all $\alpha\ge0$, and $w_\alpha\to0$
		as $\alpha\to\infty$;
		\item[(ii)] $Q$ is continuous on $[0,\infty)$, and $Q(\alpha)\to\|w_g\|_2^2$ as $\alpha\to\infty$;
		\item[(iii)] if
		\begin{equation}\label{cond:Exist-Dir}
			\|w_0-w_g\|_2<\|w_g\|_2,
		\end{equation}
		then $Q$ has a minimizer $\alpha^*\in[0,\infty)$, and $Q(\alpha^*)<\|w_g\|_2^2$.
	\end{enumerate}
	The regularizers $\frac12\|\mathbf{K}_Dw\|_2^2$ and $\mathcal{R}_{TV}$ satisfy these hypotheses, because
	$K_D$ and $K$ are injective.
\end{proposition}

\begin{proof}
	Let $\|x\|_{\mathcal{M}}^2:=x^\top\mathcal{M}_{\mathrm{LK}}x$, so that
	$D(w)-D(w_0)=\frac12\|w-w_0\|_{\mathcal{M}}^2$.

	(i) Since $R(0)=0$, we have $D(w_\alpha)+\alpha R(w_\alpha)\le D(0)$, and
	$D(0)-D(w_0)=\frac12\|w_0\|_{\mathcal{M}}^2$. Hence
	$\frac12\|w_\alpha-w_0\|_{\mathcal{M}}^2+\alpha R(w_\alpha)\le\frac12\|w_0\|_{\mathcal{M}}^2$. This gives the
	bound on $\alpha R(w_\alpha)$ and shows that $\{w_\alpha:\alpha\ge0\}$ is bounded. Let $\alpha_k\to\infty$.
	Then $R(w_{\alpha_k})\to0$, and every accumulation point $\bar w$ of $(w_{\alpha_k})$ satisfies
	$R(\bar w)=0$ by continuity of $R$, so $\bar w=0$. Hence $w_{\alpha_k}\to0$. Now let $\alpha_k\to\hat\alpha\in[0,\infty)$, and let $\bar w$ be an accumulation
	point of $(w_{\alpha_k})$. For every $w$, $D(w_{\alpha_k})+\alpha_kR(w_{\alpha_k})\le D(w)+\alpha_kR(w)$.
	Passing to the limit along the corresponding subsequence gives
	$D(\bar w)+\hat\alpha R(\bar w)\le D(w)+\hat\alpha R(w)$, so $\bar w=w_{\hat\alpha}$ by
	Lemma~\ref{lem:wellposed}. Hence $w_{\alpha_k}\to w_{\hat\alpha}$.

	(ii) This follows from (i).

	(iii) Since $Q(\alpha)\to\|w_g\|_2^2>Q(0)$, there is $A>0$ with $Q(\alpha)>Q(0)$ for all $\alpha\ge A$.
	The continuous function $Q$ attains its minimum on $[0,A]$ at some $\alpha^*$. This point minimizes
	$Q$ on $[0,\infty)$, and $Q(\alpha^*)\le Q(0)<\|w_g\|_2^2$.
\end{proof}

Condition~\eqref{cond:Exist-Dir} is sufficient but not
necessary, since a finite minimizer can exist when it fails. It cannot, however, be dropped from
Proposition~\ref{prop:exist}. Consider an abstract instance with a single pixel, $A_p=I_2$ and the
quadratic regularizer. Then $L=4I_2$ and $w_\alpha=w_0/(1+4\alpha)$. With $w_0=(1,0)^\top$ and
$w_g=(-1,0)^\top$, $Q(\alpha)=\bigl(1+(1+4\alpha)^{-1}\bigr)^2$ decreases strictly to $\|w_g\|_2^2=1$ and
never attains its infimum. For the quadratic energy, $w_\alpha$ never vanishes at a finite $\alpha$ when
$f_{\mathrm{LK}}\neq0$, because $X(\alpha)$ is invertible. The TV case behaves differently
(Section~\ref{sec:TV}).

\begin{theorem}\label{th:clg-lossE2}
	Assume (A1). Then $Q'(0)=-2\langle\mathcal{M}_{\mathrm{LK}}^{-1}Lw_0,\,w_0-w_g\rangle$. If
	\begin{equation}\label{cond:CLG-NZ}
		\bigl\langle\mathcal{M}_{\mathrm{LK}}^{-1}Lw_0,\;w_0-w_g\bigr\rangle>0,
	\end{equation}
	then $Q(\alpha)<Q(0)$ for all sufficiently small $\alpha>0$, and every minimizer of $Q$ on $[0,\infty)$
	satisfies $\alpha^*>0$.
\end{theorem}

\begin{proof}
	At $\alpha=0$, \eqref{eq:w_derivative} gives $w'_0=-\mathcal{M}_{\mathrm{LK}}^{-1}Lw_0$, and
	\eqref{eq:Q_gradient} then gives the formula for $Q'(0)$. Under~\eqref{cond:CLG-NZ}, $Q'(0)<0$, so $Q$
	decreases on a right neighborhood of $0$ and $\alpha=0$ is not a minimizer.
\end{proof}

Theorem~\ref{th:clg-lossE2} is the instance, for the CLG data term and the quadratic regularizer, of a
general sufficient condition for positivity due to Ehrhardt, Gazzola and
Scott~\cite[Corollary~1]{ehrhardt2024optimal}, which covers convex, continuously differentiable
regularizers. For an identity forward operator and a quadratic regularizer, their condition is the one
of Kunisch and Pock~\cite[Proposition~3.1]{kunisch2013bilevel}.

\begin{remark}\label{rem:first_order}
	Condition~\eqref{cond:CLG-NZ} has a first-order meaning. Since $w_\alpha=w_0+\alpha w'_0+O(\alpha^2)$ with
	$w'_0=-\mathcal{M}_{\mathrm{LK}}^{-1}Lw_0$,
	\[
		Q(\alpha)=Q(0)+2\alpha\,\langle w'_0,\,w_0-w_g\rangle+O(\alpha^2).
	\]
	Condition~\eqref{cond:CLG-NZ} therefore states that the first-order smoothing correction $w'_0$ has a
	negative component along the residual $w_0-w_g$. At each pixel,
	$w'_0(p)=-(A_p^\top A_p)^{-1}(Lw_0)(p)$, so the correction can be large where $A_p^\top A_p$ is nearly
	singular. This condition concerns the sum over all
	pixels, not each pixel separately.
\end{remark}

\begin{remark}\label{remark:neumann}
	Suppose the regularizer uses Neumann-type differences (forward differences set to zero in the last
	column and the last row), in either the quadratic or the TV case. Its kernel is then the
	space $\mathcal{Z}:=\operatorname{range}Z$ of constant flows, where
	$Z:=\begin{bmatrix}\mathbf{1}&0\\0&\mathbf{1}\end{bmatrix}\in\mathbb{R}^{2N\times2}$ and
	$\mathbf{1}\in\mathbb{R}^N$ is the vector of ones. Lemma~\ref{lem:wellposed} still applies. As
	$\alpha\to\infty$, $w_\alpha$ converges to the minimizer of $D$ on $\mathcal{Z}$,
	\[
		w_\alpha\to P_{\mathcal{M}}w_0,\qquad
		P_{\mathcal{M}}:=Z\bigl(Z^\top\mathcal{M}_{\mathrm{LK}}Z\bigr)^{-1}Z^\top\mathcal{M}_{\mathrm{LK}},
	\]
	which is the constant flow $\bigl(\sum_pA_p^\top A_p\bigr)^{-1}\sum_pA_p^\top b_p$, the least-squares
	solution of all local systems together. The proof is the argument of Proposition~\ref{prop:exist},
	with $0$ replaced by this minimizer. A sufficient condition for a finite minimizer is therefore
	\begin{equation}\label{cond:CLG-Exist-Neumann}
		\|w_0-w_g\|_2<\|P_{\mathcal{M}}w_0-w_g\|_2 .
	\end{equation}
	The $2\times2$ matrix $Z^\top\mathcal{M}_{\mathrm{LK}}Z=\sum_pA_p^\top A_p$ couples the two flow components.
	The operator $P_{\mathcal{M}}$ is the $\mathcal{M}_{\mathrm{LK}}$-orthogonal projector onto $\mathcal{Z}$, and it
	coincides with the Euclidean projector (componentwise averaging) if and only if all matrices
	$A_p^\top A_p$ are equal.
\end{remark}

\begin{theorem}\label{th:sosc}
	Assume (A1), and let $\alpha^*\ge0$ satisfy the first-order necessary conditions for
	problem~\eqref{eq:SL_CLG},
	\begin{equation}\label{eq:foc}
		\begin{cases}
			Q'(0)\ge0 & \text{if } \alpha^*=0,\\
			Q'(\alpha^*)=0 & \text{if } \alpha^*>0.
		\end{cases}
	\end{equation}
	Then $\alpha^*$ is a strict local minimizer of $Q$ on $[0,\infty)$ if either:
	\begin{enumerate}
		\item[(a)] $\alpha^*=0$ and $Q'(0)>0$; or
		\item[(b)] $\alpha^*>0$, $Lw_{\alpha^*}\neq0$, and
		\begin{equation}\label{eq:residual_bound}
			\|w_{\alpha^*}-w_g\|_2<\frac{\|w'_{\alpha^*}\|_2^2}{2\,\|X(\alpha^*)Lw'_{\alpha^*}\|_2}.
		\end{equation}
	\end{enumerate}
\end{theorem}

\begin{proof}

	(a) If $\alpha^*=0$ and $Q'(0)>0$, then $Q$ is increasing on a right neighborhood of $0$.

	(b) By~\eqref{eq:Q_gradient}, $Q''(\alpha^*)=H_1+H_2$ with
	\[
		H_1:=2\|w'_{\alpha^*}\|_2^2,\qquad H_2:=-4\langle X(\alpha^*)Lw'_{\alpha^*},\,w_{\alpha^*}-w_g\rangle .
	\]
	Since $X(\alpha^*)$ and $L$ are invertible, $Lw_{\alpha^*}\neq0$ gives $w'_{\alpha^*}\neq0$ and
	$X(\alpha^*)Lw'_{\alpha^*}\neq0$. By the Cauchy--Schwarz inequality,
	$|H_2|\le4\|X(\alpha^*)Lw'_{\alpha^*}\|_2\,\|w_{\alpha^*}-w_g\|_2$, which is smaller than $H_1$
	under~\eqref{eq:residual_bound}. Hence $Q''(\alpha^*)>0$, and together with $Q'(\alpha^*)=0$ this makes
	$\alpha^*$ a strict local minimizer.
\end{proof}

\begin{remark}\label{rem:sosc}
	The bound~\eqref{eq:residual_bound} is sufficient but not necessary: if
	$H_2\ge0$, then $Q''(\alpha^*)\ge H_1>0$ without it. In the examples of Section~\ref{sec:examples}, strict
	local optimality is checked by evaluating $Q''(\alpha^*)$ directly.
\end{remark}

Proposition~\ref{prop:exist} and Theorems~\ref{th:clg-lossE2} and~\ref{th:sosc} do not determine whether
$Q$ has a unique global minimizer (see Section~\ref{subsec:limitations}). Because $w_\alpha$ depends
nonlinearly on $\alpha$, $Q$ need not be convex, and multiple local minima cannot be excluded a priori.

\section{Total Variation Regularization}\label{sec:TV}

This section treats the TV energy~\eqref{eq:E_TV} with the gradient $\nabla=K$ of Section~\ref{sec:setting}. The TV term is not differentiable, so the
lower-level optimality condition is an inclusion rather than a linear system, and no closed form such
as~\eqref{eq:w_alpha_global} is available. 

Write $\mathcal{R}_{TV}=h\circ\mathbf{K}$ with $h(y):=\sum_p\bigl(\|y_u(p)\|_2+\|y_v(p)\|_2\bigr)$ for
$y=(y_u,y_v)\in\mathbb{R}^{4N}$, where $y_u(p),y_v(p)\in\mathbb{R}^2$ are the entries at pixel $p$. Since $h$
is finite and convex, $\partial\mathcal{R}_{TV}(w)=\mathbf{K}^\top\partial h(\mathbf{K}w)$. Hence $w_\alpha$
minimizes $E(\cdot,\alpha)$ in~\eqref{eq:E_TV} if and only if there is a dual field
$z_\alpha=(z_{\alpha,u},z_{\alpha,v})\in\mathbb{R}^{4N}$ such that, for every pixel $p$,
\begin{gather}
	A_p^\top\bigl(A_pw_\alpha(p)-b_p\bigr)=\alpha\,(\mathbf{div}\,z_\alpha)(p),\qquad
	(\mathbf{div}\,z_\alpha)(p)=\bigl((\operatorname{div}z_{\alpha,u})(p),\,(\operatorname{div}z_{\alpha,v})(p)\bigr)^\top,\label{eq:TV_oc_p}\\
	z_{\alpha,u}(p)\in\partial\|\cdot\|_2\bigl(\nabla u_\alpha(p)\bigr),\qquad
	z_{\alpha,v}(p)\in\partial\|\cdot\|_2\bigl(\nabla v_\alpha(p)\bigr).\label{eq:TV_sub}
\end{gather}
Here $\partial\|\cdot\|_2(g)=\{g/\|g\|_2\}$ if $g\neq0$, and $\partial\|\cdot\|_2(0)$ is the closed unit ball of
$\mathbb{R}^2$. In global form,
\begin{equation}\label{LL_OC_TV}
	\mathcal{M}_{\mathrm{LK}}w_\alpha-\alpha\,\mathbf{div}\,z_\alpha=f_{\mathrm{LK}},\qquad z_\alpha\in\partial h(\mathbf{K}w_\alpha).
\end{equation}
The upper-level problem is
\begin{equation}\label{eq:bilevel_CLG-TV}
	\min_{\alpha\ge0}\ Q(\alpha)=\|w_\alpha-w_g\|_2^2
	\quad\text{subject to}\quad
	w_\alpha=\operatorname*{arg\,min}_{w\in\mathbb{R}^{2N}}E(w,\alpha),
\end{equation}
and its constraint holds if and only if~\eqref{LL_OC_TV} holds for some $z_\alpha$. The dual field
$z_\alpha$ need not be unique, while $w_\alpha$ is.

\begin{corollary}\label{thm:tv_exist}
	Assume (A1), and let $\bar\alpha:=\|f_{\mathrm{LK}}\|_2/\sigma_{\min}(K)$, where $\sigma_{\min}(K)>0$ is the
	smallest singular value of $K$. Then $w_\alpha=0$ for every $\alpha\ge\bar\alpha$. Consequently
	$Q(\alpha)=\|w_g\|_2^2$ for $\alpha\ge\bar\alpha$, and problem~\eqref{eq:bilevel_CLG-TV} has a solution
	$\alpha^*\in[0,\bar\alpha]$. If in addition
	\begin{equation}\label{cond:CLG-TV-Exist}
		\|w_0-w_g\|_2<\|w_g\|_2,
	\end{equation}
	then $Q(\alpha^*)<\|w_g\|_2^2$, $\alpha^*<\bar\alpha$ and $w_{\alpha^*}\neq0$. The threshold $\bar\alpha$ is
	sufficient and need not be the smallest $\alpha$ with $w_\alpha=0$.
\end{corollary}

\begin{proof}
	Let $B:=\partial h(0)=\{y:\ \|y_u(p)\|_2\le1\text{ and }\|y_v(p)\|_2\le1\text{ for all }p\}$. Since $D$ is
	differentiable with $\nabla D(0)=-f_{\mathrm{LK}}$, the point $0$ minimizes $E(\cdot,\alpha)$ if and only if
	$f_{\mathrm{LK}}\in\alpha\,\mathbf{K}^\top B$. If $f_{\mathrm{LK}}=0$, then $w_\alpha=0$ for all $\alpha\ge0$.
	Otherwise let $\alpha\ge\bar\alpha>0$, and put $y:=\alpha^{-1}\mathbf{K}(\mathbf{K}^\top\mathbf{K})^{-1}f_{\mathrm{LK}}$,
	which is well defined because $\mathbf{K}$ is injective. Then $\mathbf{K}^\top y=f_{\mathrm{LK}}/\alpha$. By a
	singular value decomposition, $\|\mathbf{K}(\mathbf{K}^\top\mathbf{K})^{-1}\|=1/\sigma_{\min}(\mathbf{K})=1/\sigma_{\min}(K)$,
	so $\|y\|_2\le\|f_{\mathrm{LK}}\|_2/(\alpha\,\sigma_{\min}(K))\le1$. Hence $y\in B$, and $w_\alpha=0$ by uniqueness.

	By Proposition~\ref{prop:exist}(ii), $Q$ is continuous on the compact interval $[0,\bar\alpha]$, so it
	attains its minimum there at some $\alpha^*$. Since $Q(\alpha)=\|w_g\|_2^2=Q(\bar\alpha)\ge Q(\alpha^*)$ for
	$\alpha\ge\bar\alpha$, the point $\alpha^*$ solves~\eqref{eq:bilevel_CLG-TV}. Under~\eqref{cond:CLG-TV-Exist},
	$Q(\alpha^*)\le Q(0)<\|w_g\|_2^2$, which excludes both $\alpha^*\ge\bar\alpha$ and $w_{\alpha^*}=0$.
\end{proof}

In the TV case, condition~\eqref{cond:CLG-TV-Exist} is therefore not needed for existence. It guarantees
that the optimal parameter improves on the zero flow, which is the TV solution for all large $\alpha$.
Example~\ref{ex:tv_50x50} shows the difference: there~\eqref{cond:CLG-TV-Exist} fails, yet an optimal
parameter exists and improves on both the unregularized estimate and the zero flow.

\begin{remark}\label{rem:tv_neumann}
	With Neumann-type differences, $\mathcal{R}_{TV}$ vanishes on all constant flows. As $\alpha\to\infty$,
	$w_\alpha$ then converges to the constant flow $P_{\mathcal{M}}w_0$ of Remark~\ref{remark:neumann}, and
	condition~\eqref{cond:CLG-Exist-Neumann} gives a finite optimal parameter.
\end{remark}

We now give sufficient conditions under which $\alpha=0$ is not optimal, that is, under which some
$\alpha>0$ gives $Q(\alpha)<Q(0)$, so that TV regularization strictly improves on the unregularized
Lucas--Kanade estimate. The grid is assumed to have at least two rows and two columns.

For each pixel $p$, let $M_p:=(A_p^\top A_p)^{-1}$, a symmetric positive definite $2\times2$ matrix with
entries $m_{ij}(p)$, $i,j\in\{u,v\}$, where $u$ indexes the first row and column. Let
$r(p)=(r_u(p),r_v(p))^\top:=w_0(p)-w_g(p)$ be the residual of the unregularized estimate. The spatial variation of $M_p$ is measured by
$\partial_dM_p:=M_{p+e_d}-M_p$ for $d\in\{x,y\}$ when $p+e_d\in\Omega$, and by the backward difference
$\partial_dM_p:=M_p-M_{p-e_d}$ when $p$ lies in the last column ($d=x$) or the last row ($d=y$). With the
spectral norm $\|\cdot\|_{\mathrm{op}}$, let
\[
	\rho(p):=\bigl(\|\partial_xM_p\|_{\mathrm{op}}^2+\|\partial_yM_p\|_{\mathrm{op}}^2\bigr)^{1/2}.
\]
The pixels are partitioned according to which components of $w_0$ have a nonzero gradient:
\[
	\Omega_{11}=\{p:\nabla u_0(p)\ne0,\ \nabla v_0(p)\ne0\},\qquad
	\Omega_{10}=\{p:\nabla u_0(p)\ne0,\ \nabla v_0(p)=0\},
\]
\[
	\Omega_{01}=\{p:\nabla u_0(p)=0,\ \nabla v_0(p)\ne0\},\qquad
	\Omega_{00}=\{p:\nabla u_0(p)=0,\ \nabla v_0(p)=0\}.
\]
The data-alignment sums, the variation bound and the mixed-pixel bound are
\begin{align*}
	S_u&:=\sum_{p:\,\nabla u_0(p)\ne0}\Bigl[m_{11}(p)\Bigl\langle\nabla r_u(p),\tfrac{\nabla u_0(p)}{\|\nabla u_0(p)\|_2}\Bigr\rangle
	+m_{12}(p)\Bigl\langle\nabla r_v(p),\tfrac{\nabla u_0(p)}{\|\nabla u_0(p)\|_2}\Bigr\rangle\Bigr],\\
	S_v&:=\sum_{p:\,\nabla v_0(p)\ne0}\Bigl[m_{12}(p)\Bigl\langle\nabla r_u(p),\tfrac{\nabla v_0(p)}{\|\nabla v_0(p)\|_2}\Bigr\rangle
	+m_{22}(p)\Bigl\langle\nabla r_v(p),\tfrac{\nabla v_0(p)}{\|\nabla v_0(p)\|_2}\Bigr\rangle\Bigr],\\
	\mathcal{R}_{\mathrm{var}}&:=\sum_{p\in\Omega}\rho(p)\bigl(|r_u(p)|+|r_v(p)|+\|\nabla r_u(p)\|_2+\|\nabla r_v(p)\|_2\bigr),\\
	\mathcal{R}_{\mathrm{mix}}&:=\sum_{p\in\Omega_{01}}|m_{12}(p)|\,\|\nabla r_v(p)\|_2+\sum_{p\in\Omega_{10}}|m_{12}(p)|\,\|\nabla r_u(p)\|_2 .
\end{align*}
 Finally, let
\[
	D^+Q(0):=\limsup_{\alpha\downarrow0}\frac{Q(\alpha)-Q(0)}{\alpha}
\]
be the upper right Dini derivative of $Q$ at $0$. 

For a vector field $\Phi=(\Phi_u,\Phi_v)$ on $\Omega$ we write
$\partial_d\Phi(p):=\bigl((\partial_d\Phi_u)(p),(\partial_d\Phi_v)(p)\bigr)^\top$, where $\partial_x$ and
$\partial_y$ are the two components of $K$, and $\bar\Phi$ denotes the extension of $\Phi$ by zero outside
$\Omega$.

\begin{lemma}\label{lem:prod}
	Let $\Phi(p):=M_pr(p)$ for $p\in\Omega$. For every $p\in\Omega$ and $d\in\{x,y\}$,
	\[
		\partial_d\Phi(p)=M_p\,\partial_dr(p)+(\partial_dM_p)\,\bar r(p+e_d).
	\]
\end{lemma}

\begin{proof}
	If $p+e_d\in\Omega$, then
	\[\partial_d\Phi(p)=M_{p+e_d}r(p+e_d)-M_pr(p)=M_p\bigl(r(p+e_d)-r(p)\bigr)+(M_{p+e_d}-M_p)\,r(p+e_d).\]
	If $p+e_d\notin\Omega$, then $\partial_d\Phi(p)=0-M_pr(p)=M_p\,\partial_dr(p)$, and the second term vanishes
	because $\bar r(p+e_d)=0$.
\end{proof}

\begin{lemma}\label{lem:var}
	For each $p\in\Omega$ let $\zeta_u(p),\zeta_v(p)\in\mathbb{R}^2$ satisfy $\|\zeta_u(p)\|_2\le1$ and
	$\|\zeta_v(p)\|_2\le1$, and set $\eta_d(p):=\bigl([\zeta_u(p)]_d,[\zeta_v(p)]_d\bigr)^\top\in\mathbb{R}^2$ for
	$d\in\{x,y\}$. Then
	\[
		\Bigl|\sum_{p\in\Omega}\sum_{d\in\{x,y\}}\bigl\langle(\partial_dM_p)\,\bar r(p+e_d),\,\eta_d(p)\bigr\rangle\Bigr|
		\le\sqrt2\,\mathcal{R}_{\mathrm{var}} .
	\]
\end{lemma}

\begin{proof}
	Fix $p$. Since $\|\eta_x(p)\|_2^2+\|\eta_y(p)\|_2^2=\|\zeta_u(p)\|_2^2+\|\zeta_v(p)\|_2^2\le2$, the Cauchy--Schwarz
	inequality in $\mathbb{R}^2$ and then in the index $d$ gives
	\begin{align*}
		\Bigl|\sum_d\bigl\langle(\partial_dM_p)\bar r(p+e_d),\eta_d(p)\bigr\rangle\Bigr|
		&\le\sum_d\|\partial_dM_p\|_{\mathrm{op}}\|\bar r(p+e_d)\|_2\|\eta_d(p)\|_2\\
		&\le\sqrt2\,\rho(p)\max_d\|\bar r(p+e_d)\|_2 .
	\end{align*}
	In both cases of the definition, $\bar r(p+e_d)=r(p)+\partial_dr(p)$, so
	\begin{align*}
		\|\bar r(p+e_d)\|_2&\le|r_u(p)|+|r_v(p)|+|\partial_dr_u(p)|+|\partial_dr_v(p)|\\
		&\le|r_u(p)|+|r_v(p)|+\|\nabla r_u(p)\|_2+\|\nabla r_v(p)\|_2 .
	\end{align*}
	Summing over $p$ gives the claim.
\end{proof}

\begin{theorem}\label{thm:tv_nonzero}
	Assume (A1) and
	\begin{subequations}\label{cond:CLG-TV-NZ1}
		\begin{gather}
			\nabla u_0(p)=0\;\Longrightarrow\;\nabla u_g(p)=0,\quad
			\nabla v_0(p)=0\;\Longrightarrow\;\nabla v_g(p)=0
			\quad\forall p\in\Omega,
			\tag{\ref{cond:CLG-TV-NZ1}a}\label{cond:CLG-TV-NZ1a}\\
			S_u+S_v-\sqrt2\,\mathcal{R}_{\mathrm{var}}-\mathcal{R}_{\mathrm{mix}}>0.
			\tag{\ref{cond:CLG-TV-NZ1}b}\label{cond:CLG-TV-NZ1b}
		\end{gather}
	\end{subequations}
	Then $D^+Q(0)\le-2\bigl(S_u+S_v-\sqrt2\,\mathcal{R}_{\mathrm{var}}-\mathcal{R}_{\mathrm{mix}}\bigr)<0$.
	Consequently $Q(\alpha)<Q(0)$ for all sufficiently small $\alpha>0$, and $\alpha=0$ is not a solution
	of~\eqref{eq:bilevel_CLG-TV}. Every solution, which exists by Corollary~\ref{thm:tv_exist}, satisfies
	$\alpha^*>0$.
\end{theorem}

\begin{proof}
	Subtracting~\eqref{eq:TV_oc_p} at $\alpha=0$, where $A_p^\top(A_pw_0(p)-b_p)=0$, gives
	\[
		w_\alpha(p)-w_0(p)=\alpha M_p(\mathbf{div}\,z_\alpha)(p).
	\]
	With $\Phi_\alpha(p):=M_p\bigl(w_\alpha(p)+w_0(p)-2w_g(p)\bigr)$ and the symmetry of $M_p$,
	\[
		\frac{Q(\alpha)-Q(0)}{\alpha}
		=\frac1\alpha\sum_p\bigl\langle w_\alpha(p)-w_0(p),\,w_\alpha(p)+w_0(p)-2w_g(p)\bigr\rangle
		=\sum_p\bigl\langle(\mathbf{div}\,z_\alpha)(p),\Phi_\alpha(p)\bigr\rangle ,
	\]
	and the adjoint identity~\eqref{eq:adjoint} for each component yields
	\begin{equation}\label{eq:tv_quot}
		\frac{Q(\alpha)-Q(0)}{\alpha}
		=-\sum_p\Bigl[\bigl\langle\nabla\Phi_{\alpha,u}(p),z_{\alpha,u}(p)\bigr\rangle+\bigl\langle\nabla\Phi_{\alpha,v}(p),z_{\alpha,v}(p)\bigr\rangle\Bigr].
	\end{equation}

	By Proposition~\ref{prop:exist}(i), $w_\alpha\to w_0$ as $\alpha\downarrow0$, and $w_\alpha$ is bounded for
	$\alpha\in(0,1]$. Since the dual fields take values in the unit ball, the right-hand side
	of~\eqref{eq:tv_quot} is bounded, and $D^+Q(0)$ is finite. Choose $\alpha_n\downarrow0$ such that the
	quotients converge to $D^+Q(0)$, and pass to a subsequence along which $z_{\alpha_n,u}\to z_{0,u}$ and
	$z_{\alpha_n,v}\to z_{0,v}$. Since $\nabla u_{\alpha_n}(p)\to\nabla u_0(p)$ and the graph of $\partial\|\cdot\|_2$
	is closed, $z_{0,u}(p)\in\partial\|\cdot\|_2(\nabla u_0(p))$, and likewise for $v$. Moreover
	$\Phi_{\alpha_n}\to2\Phi$ with $\Phi(p):=M_pr(p)$. With
	$\eta_d(p):=\bigl([z_{0,u}(p)]_d,[z_{0,v}(p)]_d\bigr)^\top$, the limit in~\eqref{eq:tv_quot} is
	\[
		D^+Q(0)=-2\sum_p\Bigl[\langle\nabla\Phi_u(p),z_{0,u}(p)\rangle+\langle\nabla\Phi_v(p),z_{0,v}(p)\rangle\Bigr]
		=-2\sum_p\sum_d\bigl\langle\partial_d\Phi(p),\eta_d(p)\bigr\rangle .
	\]

	By Lemma~\ref{lem:prod}, the last double sum equals $T_1+T_2$ with
	\[
		T_1:=\sum_p\sum_d\bigl\langle M_p\,\partial_dr(p),\eta_d(p)\bigr\rangle,\qquad
		T_2:=\sum_p\sum_d\bigl\langle(\partial_dM_p)\,\bar r(p+e_d),\eta_d(p)\bigr\rangle ,
	\]
	and Lemma~\ref{lem:var} with $\zeta_u=z_{0,u}$ and $\zeta_v=z_{0,v}$ gives
	$|T_2|\le\sqrt2\,\mathcal{R}_{\mathrm{var}}$. The expansion of the matrix product in $T_1$ gives
	\[
		T_1=\sum_p\sum_{i\in\{u,v\}}\Bigl\langle\sum_{j\in\{u,v\}}m_{ij}(p)\,\nabla r_j(p),\,z_{0,i}(p)\Bigr\rangle .
	\]
	On $\Omega_{11}$ both subgradients are determined:
	\[z_{0,u}(p)=\nabla u_0(p)/\|\nabla u_0(p)\|_2, \quad 	z_{0,v}(p)=\nabla v_0(p)/\|\nabla v_0(p)\|_2,\]
	and the pixel contributes exactly its summands of $S_u$ and
	$S_v$. On $\Omega_{10}$ the part with $i=u$ is the $S_u$ summand. In the part with $i=v$,
	condition~\eqref{cond:CLG-TV-NZ1a} and the linearity of $\nabla$ give
	$\nabla r_v(p)=\nabla v_0(p)-\nabla v_g(p)=0$, so this part equals
	$m_{12}(p)\langle\nabla r_u(p),z_{0,v}(p)\rangle\ge-|m_{12}(p)|\,\|\nabla r_u(p)\|_2$. The set $\Omega_{01}$ is
	symmetric, with the roles of $u$ and $v$ exchanged. On $\Omega_{00}$, condition~\eqref{cond:CLG-TV-NZ1a}
	gives $\nabla r_u(p)=\nabla r_v(p)=0$, so the contribution is $0$. Hence
	$T_1\ge S_u+S_v-\mathcal{R}_{\mathrm{mix}}$.

	Consequently
	$D^+Q(0)=-2(T_1+T_2)\le-2\bigl(S_u+S_v-\mathcal{R}_{\mathrm{mix}}-\sqrt2\,\mathcal{R}_{\mathrm{var}}\bigr)$, which is
	negative by~\eqref{cond:CLG-TV-NZ1b}. Since the $\limsup$ is negative, there is $\delta>0$ with
	$\bigl(Q(\alpha)-Q(0)\bigr)/\alpha<0$ for every $\alpha\in(0,\delta)$.
\end{proof}

\begin{remark}\label{rem:dini}
	In the TV case the dual fields $z_\alpha$ need not be unique. The proof does not assume that
	$\alpha\mapsto w_\alpha$ is differentiable at $\alpha=0$. The upper Dini derivative is always finite here
	(see the proof of Theorem~\ref{thm:tv_nonzero}), and its negativity is all the argument needs.
\end{remark}

\begin{remark}\label{rem:tv_conditions}
	Condition~\eqref{cond:CLG-TV-NZ1a} is an idealized structural assumption. It asks that the ground truth
	be flat wherever a component of $w_0$ is flat. At such pixels the residual gradient of that component
	vanishes, so the undetermined subgradient does not enter the corresponding same-component term. Its
	remaining contributions, through the off-diagonal entry $m_{12}(p)$ and through the variation of
	$M_p$, are bounded by $\mathcal{R}_{\mathrm{mix}}$ and $\sqrt2\,\mathcal{R}_{\mathrm{var}}$. The condition is not expected to hold for real data in which $w_0$ has flat regions caused by the
	aperture problem. When $w_0$ is computed from noisy images, exact zeros of $\nabla u_0$ or
	$\nabla v_0$ are not expected, and~\eqref{cond:CLG-TV-NZ1a} then holds trivially. This is the case in
	Example~\ref{ex:tv_50x50}, where all pixels lie in $\Omega_{11}$.

	Condition~\eqref{cond:CLG-TV-NZ1b} compares the aggregate alignment $S_u+S_v$ between the weighted
	error gradients and the flow gradients of $w_0$ with two worst-case bounds: $\sqrt2\,\mathcal{R}_{\mathrm{var}}$
	for the spatial variation of $M_p$, and $\mathcal{R}_{\mathrm{mix}}$ for the pixels where exactly one
	component of $\nabla w_0$ vanishes. Both bounds replace the unknown dual
	vectors by their largest possible values, so the condition is conservative. It can be evaluated from
	$w_0$, the ground truth $w_g$ and the matrices $M_p$.
\end{remark}

\begin{corollary}\label{cor:tv_simple}
	Under the assumptions of Theorem~\ref{thm:tv_nonzero}, suppose in addition that
	$M_p=(A_p^\top A_p)^{-1}\equiv M$ is spatially constant with $m_{12}=0$, and that no mixed pixels occur,
	that is, $\Omega_{10}=\Omega_{01}=\varnothing$. Then $\mathcal{R}_{\mathrm{var}}=\mathcal{R}_{\mathrm{mix}}=0$, and
	condition~\eqref{cond:CLG-TV-NZ1b} reduces to
	\[
		\sum_{p:\,\nabla u_0(p)\ne0}m_{11}\Bigl\langle\nabla(u_0-u_g)(p),\tfrac{\nabla u_0(p)}{\|\nabla u_0(p)\|_2}\Bigr\rangle
		+\sum_{p:\,\nabla v_0(p)\ne0}m_{22}\Bigl\langle\nabla(v_0-v_g)(p),\tfrac{\nabla v_0(p)}{\|\nabla v_0(p)\|_2}\Bigr\rangle>0 .
	\]
	If this inequality holds, then $\alpha=0$ is not optimal.
\end{corollary}

\begin{proof}
	If $M_p$ is constant, both the forward and the backward differences $\partial_dM_p$ vanish, so $\rho=0$
	and $\mathcal{R}_{\mathrm{var}}=0$. Without mixed pixels, $\mathcal{R}_{\mathrm{mix}}=0$. With $m_{12}=0$ the
	displayed inequality is condition~\eqref{cond:CLG-TV-NZ1b}, and the conclusion follows from
	Theorem~\ref{thm:tv_nonzero}.
\end{proof}

\begin{remark}\label{rem:weighted_conditions}
	For the weighted loss $Q_W$ of Remark~\ref{rem:weighted}, the results of Sections~\ref{sec:quad}
	and~\ref{sec:TV} hold with the following changes, and the proofs are otherwise unchanged. In
	Proposition~\ref{prop:exist}, the limit of $Q_W$ is $\|Ww_g\|_2^2$ and condition~\eqref{cond:Exist-Dir}
	becomes $\|W(w_0-w_g)\|_2<\|Ww_g\|_2$. The additional condition of Corollary~\ref{thm:tv_exist} changes in
	the same way. In Theorem~\ref{th:clg-lossE2}, condition~\eqref{cond:CLG-NZ} becomes
	$\langle\mathcal{M}_{\mathrm{LK}}^{-1}Lw_0,\,W^2(w_0-w_g)\rangle>0$. Theorem~\ref{th:sosc} applies with the first-order
	conditions~\eqref{eq:foc} written for $Q_W$. In part~(b), $Lw_{\alpha^*}\neq0$ is replaced by $Ww'_{\alpha^*}\neq0$,
	which follows from $Lw_{\alpha^*}\neq0$ when $\omega>0$ at every pixel but not for a mask,
	and~\eqref{eq:residual_bound} is replaced by
	$2\,\|WX(\alpha^*)Lw'_{\alpha^*}\|_2\,\|W(w_{\alpha^*}-w_g)\|_2<\|Ww'_{\alpha^*}\|_2^2$. For TV, let
	$\widetilde M_p:=\omega(p)^2M_p$. Theorem~\ref{thm:tv_nonzero} holds with $\widetilde M_p$ in place of $M_p$ in
	$S_u$, $S_v$, $\rho$, $\mathcal{R}_{\mathrm{var}}$ and $\mathcal{R}_{\mathrm{mix}}$, and
	Corollary~\ref{cor:tv_simple} then requires $\widetilde M_p$, not $M_p$, to be spatially constant, with
	$m_{ij}$ the entries of $\widetilde M_p$.
\end{remark}

\section{Illustrative Examples}
\label{sec:examples}

We present numerical examples that illustrate the theoretical conditions, including cases in which the sufficient conditions fail. These examples are not intended as a performance evaluation or comparison with state-of-the-art methods. Two scales are considered: a $3\times3$ grid (and a $4\times4$ grid in Example~\ref{ex:tv_designed_mix}), which allows direct inspection of all conditions, and a $50\times50$ grid, which exercises the framework at a more realistic scale with noise and subpixel motion.

In the examples built from images, frame $I_2$ is produced by bilinear backward warping of $I_1$
with a constant global translation. Gaussian noise is then added to $I_2$, and the intensities are
clipped to $[0,1]$. The spatial derivatives $I_x,I_y$ are central differences at interior pixels and one-sided differences
at the border, computed on each frame and averaged over the two frames. Each window $\Omega_p$ is the
$5\times5$ block of pixels centered at $p$, truncated at the image border. All reported numerical
values of these examples are illustrative and depend on the particular noise realization.
In the examples on small grids, the image derivatives or the matrices $A_p^\top A_p$ are
prescribed directly, and (A1) holds without modification.

The theoretical results (Proposition~\ref{prop:exist}, Theorems~\ref{th:clg-lossE2}, \ref{th:sosc} and~\ref{thm:tv_nonzero}, and Corollary~\ref{thm:tv_exist}) assume
$\mathcal{M}_{\mathrm{LK}}\succ0$, which follows from full-rank LK patches. In the $50\times50$ experiments, the
stabilization $A_p^\top A_p \mapsto A_p^\top A_p + 10^{-2}I$ ensures positive definiteness in low-texture regions
where the rank condition may not hold strictly. The theorems apply to the stabilized $\mathcal{M}_{\mathrm{LK}}$.

Quadratic regularization uses the matrix $L$ and TV regularization uses the gradient $K$ of Section~\ref{sec:setting}.

For quadratic regularization, the lower-level solution
$w_\alpha = (\mathcal{M}_{\mathrm{LK}} + \alpha L)^{-1} f_{\mathrm{LK}}$ is computed by one linear solve. For TV, the lower-level problem is solved by the Chambolle--Pock primal--dual
algorithm~\cite{chambolle2011afirst} with the exact proximal map of the data term, which is a $2\times2$
linear solve per pixel. The step sizes start from $\tau=\sigma=0.99/\sqrt8$, so that $\tau\sigma\|K\|^2<1$,
and follow the accelerated rule for a strongly convex data term. Each solve stops when the relative
primal--dual gap falls below a prescribed tolerance or at an iteration limit, which is $10^5$ for the
sweep solves. The sweep tolerances are given in each example.

In Examples~\ref{ex:quad_3x3} and~\ref{ex:quad_50x50} the optimal parameter is located by bounded
one-dimensional minimization (Brent's method), and sweeps of $Q$ over a grid of $\alpha$ values produce the
plotted curves. In Examples~\ref{ex:tv_3x3} and~\ref{ex:tv_50x50}, $Q$ is sampled on a grid of $\alpha$ values
with warm starts and then on a finer grid around the coarse minimum, and a confirmation run at $\alpha^*$
uses the tolerance $10^{-12}$. Example~\ref{ex:sosc_constructed} prescribes the parameter and checks it by a
grid sweep, and Example~\ref{ex:tv_designed_mix} evaluates the positivity condition without a parameter
search.

\subsection{Quadratic Regularization}
\label{subsec:quad_examples}

Examples~\ref{ex:quad_3x3} and~\ref{ex:quad_50x50} satisfy the existence
condition~\eqref{cond:Exist-Dir} and have finite minimizers. Remark~\ref{rem:high_noise_quad} gives a
case in which~\eqref{cond:Exist-Dir} fails and a finite minimizer still exists, and the one-pixel
instance after Proposition~\ref{prop:exist} shows that without~\eqref{cond:Exist-Dir} the infimum need
not be attained. A third, constructed example illustrates the second-order condition of
Theorem~\ref{th:sosc}(b).

\begin{example}
	\label{ex:quad_3x3}
	We consider a synthetic $3 \times 3$ image pair with per-pixel spatial gradients
	and temporal derivatives given by
	\[
		\nabla I =
		\begin{bmatrix}
			(1.0, 0.0) & (0.0, 1.0)  & (1.0, 1.0) \\
			(0.5, 0.5) & (0.7, 0.3)  & (0.0, 0.8) \\
			(0.2, 0.9) & (1.0, -0.5) & (0.4, 0.6)
		\end{bmatrix},
		\qquad
		I_t =
		\begin{bmatrix}
			-0.1 & -0.2 & -0.3 \\
			-0.4 & -0.5 & -0.6 \\
			-0.7 & -0.8 & -0.9
		\end{bmatrix},
	\]
	and a spatially varying ground-truth flow field
	\[
		w_g =
		\begin{bmatrix}
			(0.1, 0.1) & (0.2, 0.1) & (0.3, 0.2) \\
			(0.2, 0.3) & (0.4, 0.3) & (0.5, 0.4) \\
			(0.3, 0.5) & (0.6, 0.5) & (0.7, 0.6)
		\end{bmatrix}.
	\]

	For this toy example, the entire $3 \times 3$ grid is used as the local
	neighborhood for each pixel, so that all blocks $A_p^\top A_p$ and vectors
	$A_p^\top b_p$ are identical. We obtain
	\[
		A_p^\top A_p =
		\begin{bmatrix}
			3.940 & 1.380 \\
			1.380 & 4.400
		\end{bmatrix},
		\qquad
		A_p^\top b_p =
		\begin{bmatrix}
			2.250 \\
			2.100
		\end{bmatrix}.
	\]
	The constant Lucas--Kanade solution is therefore $w_0(p) \approx (0.454,\,0.335)$ at every pixel.

	\medskip
	We first verify the conditions of Proposition~\ref{prop:exist} and Theorem~\ref{th:clg-lossE2}. 
	For this example,
	\[
		Q(0)=\|w_0-w_g\|_2^2 \approx 0.648
		\quad\text{and}\quad
		\|w_g\|_2^2 \approx 2.79,
	\]
	so $Q(0)<\|w_g\|_2^2$. Hence a finite optimal regularization parameter exists under
	Dirichlet boundary conditions.

	\medskip
	The inner product
	\[
		\left\langle
		\mathcal{M}_{\mathrm{LK}}^{-1}L w_0,\,
		w_0 - w_g
		\right\rangle
		\approx 0.110 > 0,
	\]
	so condition~\eqref{cond:CLG-NZ} is satisfied and $\alpha^*>0$.
	Although $w_0$ is spatially uniform, $Lw_0\neq0$. Since $w_0$ is constant and the grid spacing is
	$1$, $(Lw_0)(p)=(4-k_p)\,w_0(p)$, where $k_p$ is the number of in-grid neighbors of $p$, that is, $k_p=2$, $3$ and $4$ at the corners, the edge midpoints and the
	center. The inner product is therefore
	$\sum_p(4-k_p)\,\langle(A_p^\top A_p)^{-1}w_0(p),\,w_0(p)-w_g(p)\rangle\approx0.1102$.

	\medskip
	Numerical minimization of
	$Q(\alpha) = \|w_\alpha - w_g\|_2^2$
	gives
	\[
		\alpha^* \approx 0.585, \qquad
		Q(\alpha^*) \approx 0.597,
	\]
	which is an error reduction of approximately $7.96\%$.
	Figure~\ref{fig:Qalpha} shows $Q(\alpha)$.

\medskip
At the computed minimizer, the first-order condition is confirmed:
$|Q'(\alpha^*)| < 10^{-7}$.
Moreover, $L w_{\alpha^*} \neq 0$ and direct numerical evaluation gives
$Q''(\alpha^*) \approx 0.204 > 0$.
Hence $\alpha^*$ is a strict local minimizer, although the residual bound~\eqref{eq:residual_bound} of Theorem~\ref{th:sosc}(b) is not satisfied here ($0.772>0.521$) (see Remark~\ref{rem:sharpness_quad}).
Table~\ref{tab:ex_quad_3x3} summarizes the key numerical values.

	\begin{table}[h]
		\centering
		\caption{Key numerical values for Example~\ref{ex:quad_3x3}.}
		\label{tab:ex_quad_3x3}
		\begin{tabular}{lc}
			\toprule
			Quantity & Value \\
			\midrule
			$Q(0)<\|w_g\|_2^2$ (Dirichlet existence, \eqref{cond:Exist-Dir}) & Yes \\
			$\langle \mathcal{M}_{\mathrm{LK}}^{-1}L w_0,\, w_0 - w_g \rangle$ & $0.110$ \\
			Condition~\eqref{cond:CLG-NZ} satisfied & Yes \\
			Optimal $\alpha^*$ & $0.585$ \\
			$Q(0) = \|w_0 - w_g\|_2^2$ & $0.648$ \\
			$\|w_g\|_2^2$ & $2.79$ \\
			$Q(\alpha^*)$ & $0.597$ \\
			Error reduction & $7.96\%$ \\
			$|Q'(\alpha^*)|$ & $<10^{-7}$ \\
			$Q''(\alpha^*)$ & $0.204$ \\
			$L w_{\alpha^*} \neq 0$ & Yes \\
			Strict local minimizer ($Q''(\alpha^*)>0$, direct evaluation) & Yes \\
			\bottomrule
		\end{tabular}
	\end{table}

	\begin{figure}[!t]
		\centering
		\includegraphics[width=0.75\textwidth]{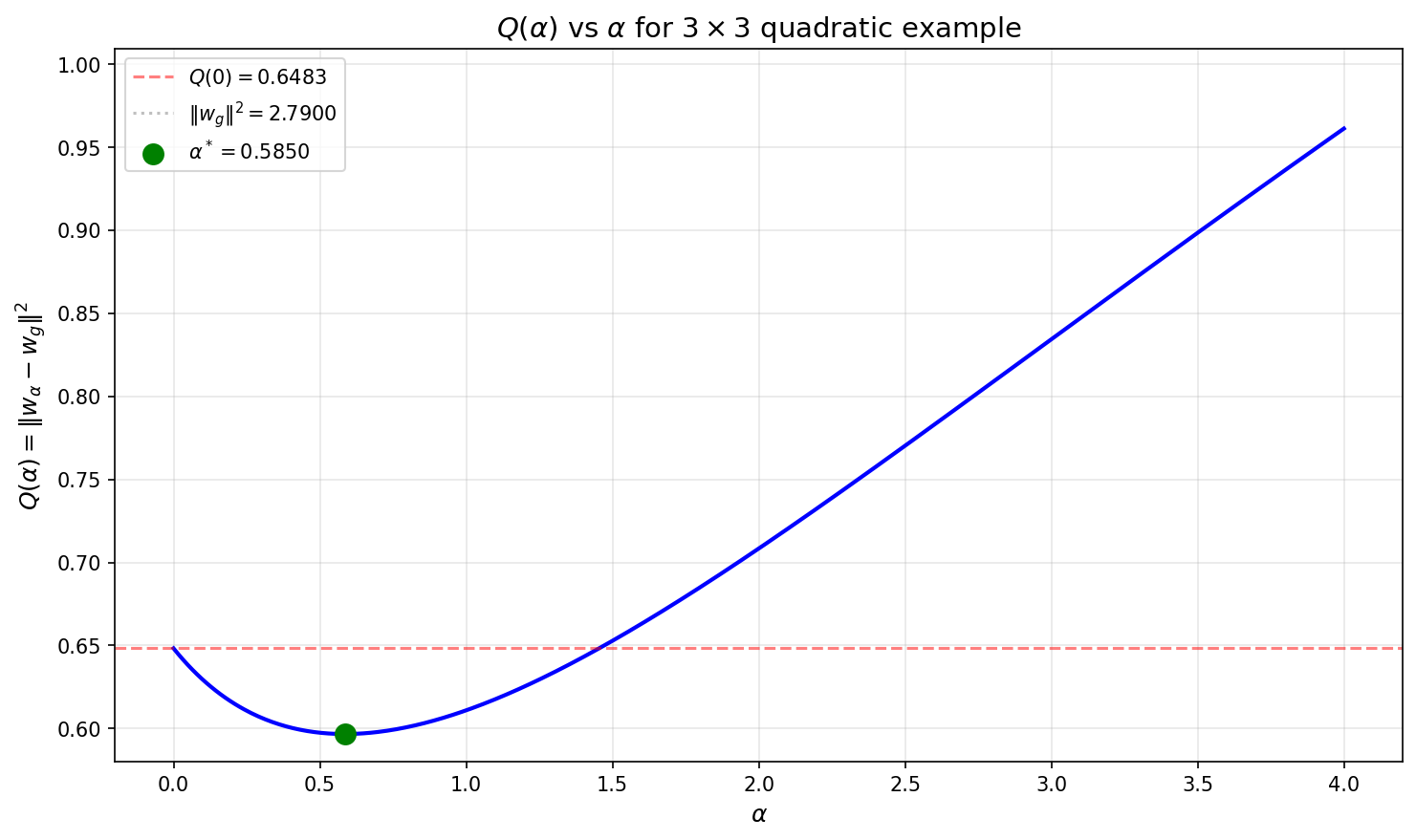}
		\caption{Objective $Q(\alpha)=\|w_\alpha-w_g\|_2^2$ versus the regularization
			parameter $\alpha$ for the $3\times3$ quadratic example.
			The curve decreases for small $\alpha>0$, which confirms the prediction
			$\alpha^*>0$ from condition~\eqref{cond:CLG-NZ}}
		\label{fig:Qalpha}
	\end{figure}
\end{example}

\begin{example}
	\label{ex:quad_50x50}
	We consider a $50\times50$ scene where frame $I_1$ is a smooth random texture obtained
	by applying a Gaussian blur (standard deviation $3$ pixels) to white noise and normalizing to $[0,1]$. Frame $I_2$ is
	the warp of $I_1$ by the constant global translation
	$w_g=(0.7,-0.3)$ at every pixel, with border values replicated, and the noise level is
	$\sigma=0.10$. Both random fields come from NumPy's
	\texttt{default\_rng} with seed $154$.
	
	Under Dirichlet boundary conditions,
	\[
		\lim_{\alpha\to\infty} Q(\alpha) \;=\; \|w_g\|_2^2 \;=:\; Q_\infty \;=\; 1450.
	\]

	\medskip
	We verify the conditions of Proposition~\ref{prop:exist} and Theorem~\ref{th:clg-lossE2}.
	For the Dirichlet finite‑existence condition~\eqref{cond:Exist-Dir},
	\[
		Q(0)=\|w_0-w_g\|_2^2 \approx 1294.8
		\;<\; \|w_g\|_2^2 = 1450.0 = Q_\infty,
	\]
	so condition~\eqref{cond:Exist-Dir} is \emph{satisfied} and
	Proposition~\ref{prop:exist} guarantees a finite optimal regularization parameter.

	\medskip
	For the positivity test (condition~\eqref{cond:CLG-NZ}),
	\[
		\left\langle
		\mathcal{M}_{\mathrm{LK}}^{-1}\,L \, w_0,\,
		w_0 - w_g
		\right\rangle
		\approx 1.24\times 10^{4} \;>\; 0,
	\]
	so condition~\eqref{cond:CLG-NZ} is satisfied and $Q'(0)<0$, which guarantees $\alpha^\ast>0$.

	\medskip
	A log‑scale sweep over $\alpha\in[10^{-2},10^{2}]$ brackets the minimizer, and the bounded search locates it at
	\[
		\alpha^* \approx 0.532,\qquad Q(\alpha^*) \approx 896.4,
	\]
	a $30.8\%$ reduction versus $Q(0)$.
	Figure~\ref{fig:ex2-dirichlet}(d)--(f) shows the curve $Q(\alpha)$ and the pointwise error maps
	$|w_0-w_g|$ and $|w(\alpha^*)-w_g|$.

	\medskip
	At the computed minimizer: $|Q'(\alpha^*)|\approx 0$,
	$Lw_{\alpha^*}\neq 0$, and
	$Q''(\alpha^*)\approx 230>0$.
	The curvature decomposes as $Q''(\alpha^*)=H_1(\alpha^*)+H_2(\alpha^*)$
	with $H_1(\alpha^*)\approx104$ and $H_2(\alpha^*)\approx126\geq0$,
	so $Q''(\alpha^*)\geq H_1(\alpha^*)>0$ holds without invoking the
	residual bound~\eqref{eq:residual_bound}.
	Six further realizations of the same generator (seeds $0$--$3$, $42$ and $99$) give $Q(0)$ between
	$1092$ and $1323$, $\alpha^*$ between $0.52$ and $0.68$, and reductions between $30.0\%$ and $38.2\%$.
	At every seed $H_2(\alpha^*)>0$ and the residual bound~\eqref{eq:residual_bound} fails. This small sample
	describes the variability of the realization and is not a robustness study.
	Hence $Q''(\alpha^*)>0$ directly confirms that $\alpha^*$ is a strict local minimizer (see Remark~\ref{rem:sharpness_quad}).
	Table~\ref{tab:ex_quad_50x50} summarizes the key values.

	\begin{table}[h]
		\centering
		\caption{Key numerical values for Example~\ref{ex:quad_50x50}.}
		\label{tab:ex_quad_50x50}
		\begin{tabular}{lc}
			\toprule
			Quantity & Value \\
			\midrule
			$Q(0)=\|w_0-w_g\|_2^2$ & $1294.8$ \\
			$\|w_g\|_2^2=Q_\infty$ & $1450.0$ \\
			$Q(0)<\|w_g\|_2^2$ (Dirichlet existence, \eqref{cond:Exist-Dir}) & Yes \\
			$\langle \mathcal{M}_{\mathrm{LK}}^{-1}L w_0,\,
			w_0 - w_g \rangle$ (positivity, \eqref{cond:CLG-NZ}) & $1.24\times 10^{4}$ \\
			$\alpha^*$ (interior minimum, guaranteed by Prop.~\ref{prop:exist} and Thm.~\ref{th:clg-lossE2}) & $0.532$ \\
			$Q(\alpha^*)$ & $896.4$ \\
			Reduction vs.\ $\alpha=0$ & $30.8\%$ \\
			$|Q'(\alpha^*)|$ & $\approx 0$ \\
			$Q''(\alpha^*)=H_1(\alpha^*)+H_2(\alpha^*)$ & $230$ \\
			$H_2(\alpha^*)=-4\langle X(\alpha^*)Lw'_{\alpha^*},\,w_{\alpha^*}-w_g\rangle$ & $126$ \\
			$H_2\geq 0 \Rightarrow Q''\geq H_1>0$ (no bound needed) & Yes \\
			Residual bound~\eqref{eq:residual_bound} satisfied & No \\
			$L w_{\alpha^*} \neq 0$ & Yes \\
			Strict local minimizer ($Q''(\alpha^*)>0$, $H_2\geq0$ route) & Yes \\
			\bottomrule
		\end{tabular}
	\end{table}

	\begin{figure}[!t]
		\centering
		\includegraphics[width=0.95\textwidth]{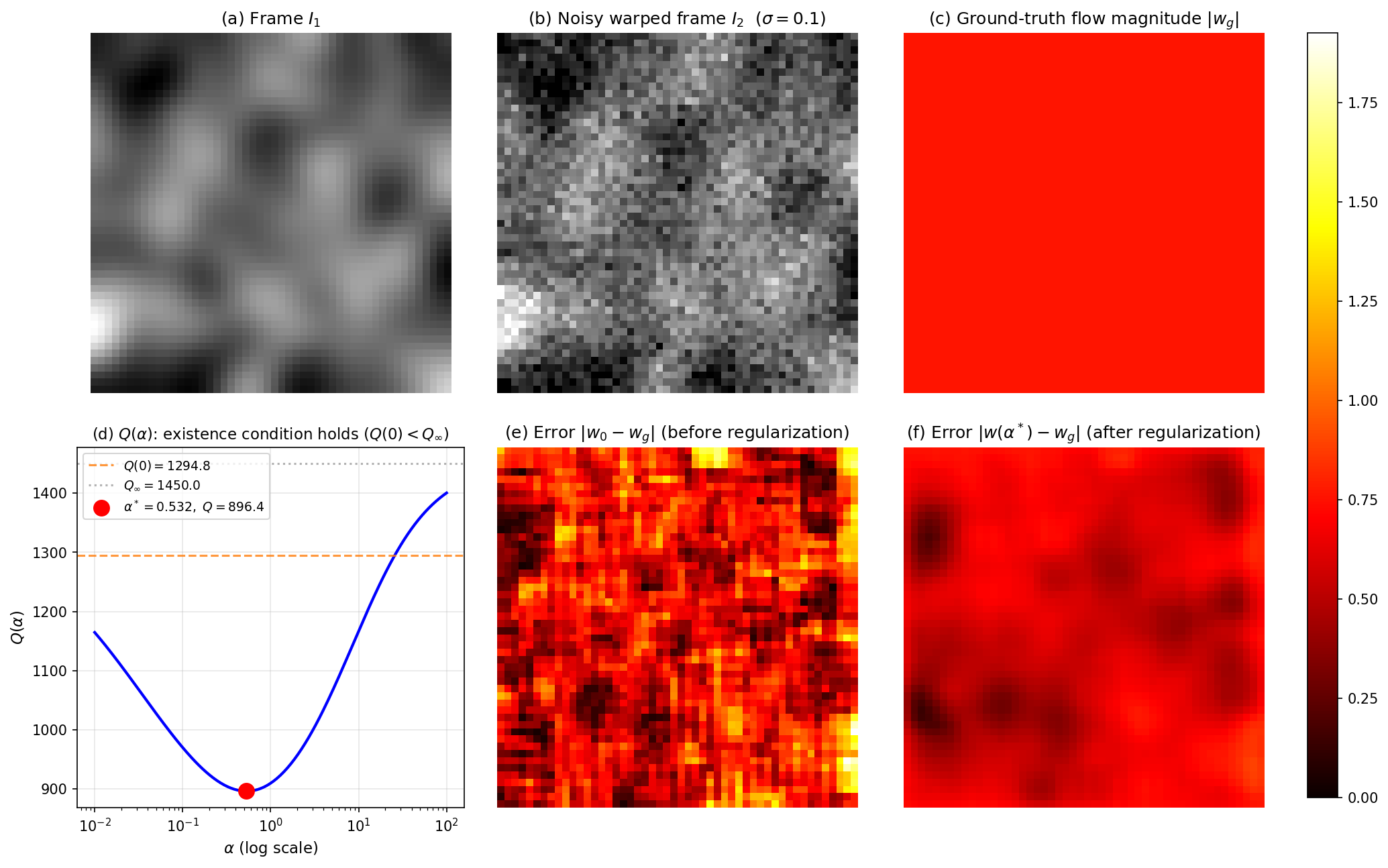}
		\caption{\textbf{Example~\ref{ex:quad_50x50}:}
			(a) Frame $I_1$ (smooth random texture).
			(b) Noisy warped frame $I_2$ ($\sigma=0.10$).
			(c) Ground‑truth flow magnitude $|w_g|$ (uniform: global translation
			$(0.7,-0.3)$, magnitude $\approx0.76$).
			(d) $Q(\alpha)$ on a log scale, with the minimizer $\alpha^*\approx0.532$ (red marker).
			(e) Pointwise error $|w_0-w_g|$ before regularization.
			(f) Pointwise error $|w(\alpha^*)-w_g|$ after regularization
			($30.8\%$ reduction in $Q$)}
		\label{fig:ex2-dirichlet}
	\end{figure}

	\medskip

	Figure~\ref{fig:noise_sensitivity} shows the noise sensitivity of the
	quadratic bilevel improvement. The improvement ratio $(Q(0)-Q(\alpha^*))/Q(0)$ rises from ${\approx}18\%$
	at $\sigma=0.02$ to ${\approx}37\%$ at $\sigma=0.10$ and then stabilizes near $36\%$.
	The figure reports medians over five further seeds ($42$--$46$) of the same generator, so its value
	at $\sigma=0.10$ differs from the $30.8\%$ of the single realization above. For this scene and these
	seeds, the benefit persists over the tested range of noise levels.

	\begin{figure}[!t]
	\centering
	\includegraphics[width=0.85\textwidth]{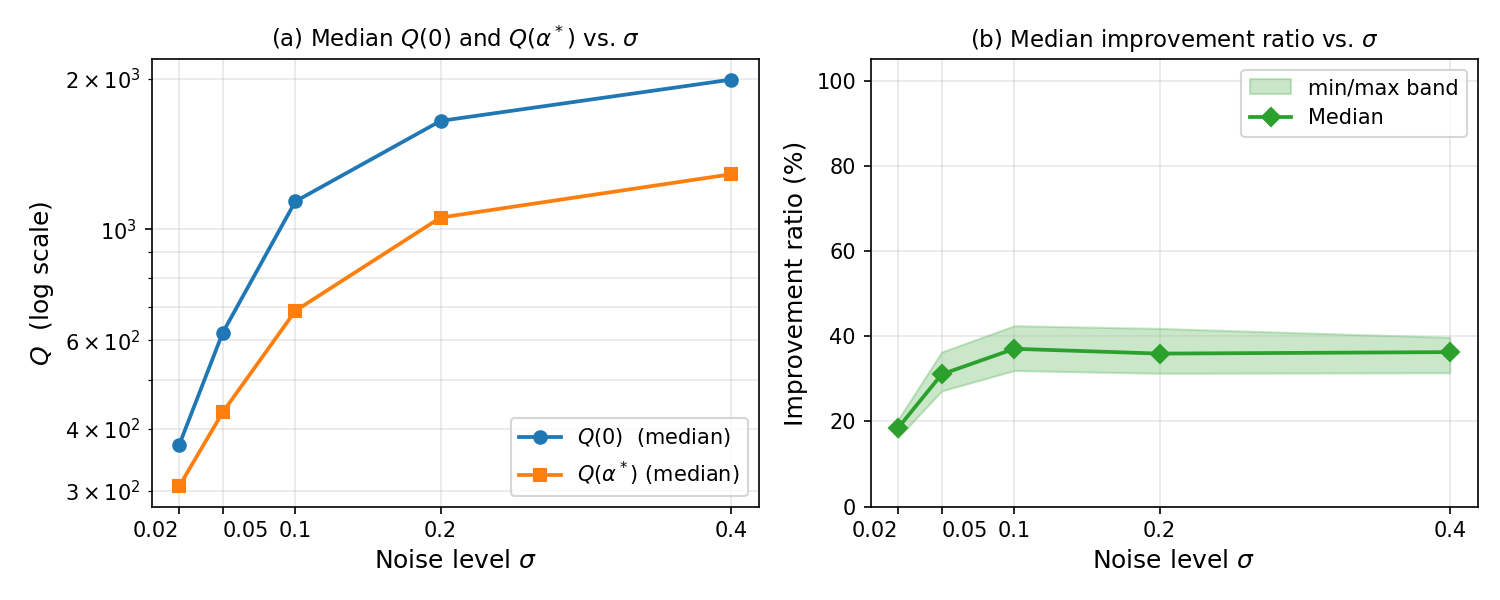}
	\caption{Noise-sensitivity of the bilevel improvement for
	Example~\ref{ex:quad_50x50} ($50\times50$ textured global-translation scene,
	quadratic regularization, five random seeds).
	(a)~Median $Q(0)$ and $Q(\alpha^*)$ versus noise level $\sigma$ (log $y$-scale).
	(b)~Median improvement ratio $(Q(0)-Q(\alpha^*))/Q(0)$ with min/max band over seeds.
	The band width remains below $11$ percentage points over the five seeds}
	\label{fig:noise_sensitivity}
	\end{figure}
\end{example}

\begin{remark}
	\label{rem:high_noise_quad}
	The existence condition~\eqref{cond:Exist-Dir} is sufficient but not necessary.
	We retain the textured scene of Example~\ref{ex:quad_50x50} (same generator and seed) and raise the
	noise to $\sigma=0.50$, at which $880$ of the $2500$ pixels of $I_2$ are clipped. This gives $Q(0)\approx 2089.0 > \|w_g\|_2^2=1450.0$, so
	condition~\eqref{cond:Exist-Dir} \emph{fails}.
	Nevertheless, numerical optimization finds a finite interior minimizer at
	$\alpha^*\approx 13.0$ with $Q(\alpha^*)\approx 1425.4 < Q_\infty=1450.0$,
	a $31.8\%$ reduction.
In the realizations computed here, condition~\eqref{cond:Exist-Dir} holds at $\sigma=0.10$ (seed $154$ and the six further seeds of
	Example~\ref{ex:quad_50x50}). It fails for all five seeds of the noise study (Figure~\ref{fig:noise_sensitivity}) at
	$\sigma=0.20$ and at $\sigma=0.40$, where $Q(0)$ lies between $1609$ and $2126$, and at $\sigma=0.50$. In all
	these cases a finite minimizer still reduces the error by $31\%$ to $42\%$. Noise inflates $Q(0)$, and in
	these realizations the condition already fails at $\sigma=0.20$, although regularization still reduces
	the error.
\end{remark}

\begin{example}
	\label{ex:sosc_constructed}
	We reuse the $3\times3$ quadratic setup of Example~\ref{ex:quad_3x3} (same
	$\mathcal{M}_{\mathrm{LK}}$, $L$, $f_{\mathrm{LK}}$) and fix
	$\alpha^*=1$. The flow $w^*$ and its sensitivity $d=w'_{\alpha^*}$
	are computed from the lower-level system, and $w_g$ is then chosen so that
	$\alpha^*$ satisfies both the first-order condition and the residual bound of
	Theorem~\ref{th:sosc}(b).

	The key quantities are defined by
	\[
	\begin{gathered}
		X^*:=(\mathcal{M}_{\mathrm{LK}}+\alpha^*L)^{-1},\qquad w^*:=X^*f_{\mathrm{LK}},\\
		d:=w'_{\alpha^*}=-X^*Lw^*,\qquad z:=X^*Ld,\\
		v:=\frac{P_{d^\perp}z}{\|P_{d^\perp}z\|_2},\qquad \varepsilon:=\frac{\|d\|_2^2}{4\|z\|_2},\qquad
		w_g:=w^*-\varepsilon v,
	\end{gathered}
	\]
	where $P_{d^\perp}z:=z-\tfrac{\langle z,d\rangle}{\|d\|_2^2}\,d$ is the component of $z$
	orthogonal to $d$. Since $Lw^*\neq0$ we have $d\neq0$. The assumption
	$P_{d^\perp}z\neq0$ ensures $v$ is well-defined, with $\|v\|_2=1$,
	$v\perp d$, and $\langle z,v\rangle>0$.

	Because $w^*-w_g=\varepsilon v$ and $v\perp d$, the first-order condition gives
	$Q'(\alpha^*)=2\langle d,\varepsilon v\rangle=0$, so $\alpha^*$ is a
	critical point (cf.~\eqref{eq:foc}).
	By construction $\varepsilon=\tfrac{\|d\|_2^2}{4\|z\|_2}<\tfrac{\|d\|_2^2}{2\|z\|_2}$,
	so the residual bound~\eqref{eq:residual_bound} holds.
	Moreover $H_2(\alpha^*)=-4\varepsilon\langle z,v\rangle<0$. The example therefore lies
	in the regime where the curvature contribution $H_2$ is negative, and
	Theorem~\ref{th:sosc}(b) guarantees that $\alpha^*$ is a strict local minimizer.

	\medskip\noindent
	For the $3\times3$ setup with $\alpha^*=1$, numerical verification gives
	$\|d\|_2\approx0.262$ and $\|z\|_2\approx0.0626$. A sweep of $Q$ over $[0,6]$ with step $10^{-3}$ attains its smallest
	value at $\alpha=1$.
	Conditions~\eqref{cond:Exist-Dir} and~\eqref{cond:CLG-NZ} both hold for this constructed $w_g$.
	Table~\ref{tab:sosc} collects the remaining quantities.

	\begin{table}[h]
		\centering
		\caption{Key numerical values for the auxiliary Example~\ref{ex:sosc_constructed}.}
		\label{tab:sosc}
		\begin{tabular}{lc}
			\toprule
			Quantity                                                            & Value                \\
			\midrule
			$Q(0)<\|w_g\|_2^2$ (Dirichlet existence, \eqref{cond:Exist-Dir})  & Yes                  \\
			Condition~\eqref{cond:CLG-NZ} satisfied                            & Yes                  \\
			Prescribed $\alpha^*$                                        & $1.000$              \\
			$Q(0)=\|w_0-w_g\|_2^2$                                             & $0.166$              \\
			$\|w_g\|_2^2$                                                      & $1.764$              \\
			$Q(\alpha^*)$                                                & $0.075$              \\
			Error reduction                                                    & $54.9\%$             \\
			$Q'(\alpha^*)$                                               & $\approx 0$          \\
			$H_1(\alpha^*)=2\|d\|_2^2$                                   & $0.137$              \\
			$H_2(\alpha^*)=-4\langle z,w^*-w_g\rangle$              & $-0.026$             \\
			$Q''(\alpha^*)=H_1+H_2$                                      & $0.111$              \\
			$\|w^*-w_g\|_2$ (residual $\varepsilon$)                     & $0.274$              \\
			$\|d\|_2^2/(2\|z\|_2)$ (SOSC bound)                               & $0.548$              \\
			Residual $<$ bound~\eqref{eq:residual_bound}                       & Yes                  \\
			Strict local minimizer ($Q''(\alpha^*)>0$, residual bound)                   & Yes                  \\
			\bottomrule
		\end{tabular}
	\end{table}
\end{example}

\begin{remark}
	\label{rem:sharpness_quad}
	Examples~\ref{ex:quad_3x3}--\ref{ex:sosc_constructed} clarify the status of the second-order condition
	in Theorem~\ref{th:sosc}. In Example~\ref{ex:quad_3x3} the residual bound~\eqref{eq:residual_bound} is not
	satisfied, yet direct evaluation confirms $Q''(\alpha^*)>0$. The bound is conservative
	at this instance.
	In Example~\ref{ex:quad_50x50} the coupling term satisfies $H_2(\alpha^*)\geq 0$
	at all seven computed realizations, so $Q''(\alpha^*)=H_1(\alpha^*)+H_2(\alpha^*)
	\geq H_1(\alpha^*)>0$ follows automatically. The residual bound, which also fails
	here, is not needed to certify strict local optimality.
	The auxiliary Example~\ref{ex:sosc_constructed} provides the complementary picture:
	the bound \emph{is} satisfied and the coupling
	term is adverse ($H_2<0$), so the residual bound certifies
	$H_1>|H_2|$ without direct evaluation of $Q''$.
\end{remark}

\subsection{Total Variation}
\label{subsec:tv_examples}

\begin{example}
	\label{ex:tv_3x3}
	We consider a $3 \times 3$ grid with a ground truth flow that exhibits a step pattern, shared by both components:
	\[
		w_g =
		\begin{bmatrix}
			(0,0) & (0,0) & (0,0) \\
			(0,0) & (1,1) & (1,1) \\
			(0,0) & (1,1) & (1,1)
		\end{bmatrix}.
	\]
	To construct an unregularized solution $w_0$ that satisfies the structural assumptions of Theorem~\ref{thm:tv_nonzero}, we introduce oscillatory perturbations around the step:
	\[
		w_0 =
		\begin{bmatrix}
			(-1,\,-1) & (0.5,\,0.5)  & (-0.5,\,-0.5) \\
			(0.5,\,0.5)  & (0.5,\,0.5)  & (1.5,\,1.5) \\
			(-0.5,\,-0.5) & (1.5,\,1.5)  & (0.5,\,0.5)
		\end{bmatrix}.
	\]
	In this example $u_0=v_0$ and $u_g=v_g$ at every pixel. With the gradient $K$, $\nabla u_0(p)\neq0$ at all nine pixels, so every pixel lies in $\Omega_{11}$ and
condition~\eqref{cond:CLG-TV-NZ1a} holds trivially.

	This example uses the specialization $\mathcal{M}_{\mathrm{LK}}=I$ (identity data term),
	so the lower-level energy reduces to
	\[
		E(w, \alpha) = \frac{1}{2} \|w - w_0\|_2^2 + \alpha \, \mathrm{TV}(w).
	\]
	With $A_p^\top A_p=I$ at every pixel, Corollary~\ref{cor:tv_simple} applies ($M_p=I$, $m_{12}=0$, no mixed pixels), and condition~\eqref{cond:CLG-TV-NZ1b} reduces to $S_u+S_v>0$.

	We first verify the conditions of Corollary~\ref{thm:tv_exist}. Numerically,
	$Q(0)=\|w_0-w_g\|_2^2=6.00$ and $\|w_g\|_2^2=8.00$.
	An optimal parameter exists by Corollary~\ref{thm:tv_exist}. Since $Q(0)<\|w_g\|_2^2$,
	condition~\eqref{cond:CLG-TV-Exist} is satisfied, so the optimum also improves on the zero flow. The
	computed flow is zero for every sampled $\alpha\ge0.72$, well below the sufficient threshold
	$\bar\alpha=\|w_0\|_2/\sigma_{\min}(K)\approx5.94$ of the corollary.

	The alignment sum $S_u+S_v \approx 17.91 > 0$ satisfies condition~\eqref{cond:CLG-TV-NZ1b},
	which by Theorem~\ref{thm:tv_nonzero} implies $D^+Q(0) < 0$ and hence $\alpha^* > 0$.
	The forward-difference estimate $D^+Q(0)\approx -35.75 < 0$ (step $10^{-4}$) confirms this result,
	consistent with $-2(S_u+S_v) \approx -35.81$ from the theoretical formula.

A sweep of $101$ uniform points on $[0,2]$, followed by $81$ points around the coarse
	minimum, gives $\alpha^* \approx 0.211$ with $Q(\alpha^*) \approx 2.558$, an error reduction of
	approximately $57.4\%$. The sweep tolerance is $10^{-10}$. The largest final relative gap was $1.1\times10^{-9}$. A confirmation run at $\alpha^*$ agrees to five digits.
	The behavior of $Q(\alpha)$ is shown in Figure~\ref{fig:tv_qalpha_ex}.
	TV regularization reduces the oscillations of $w_0$. On this small grid it recovers the step of
	$w_g$ only partly, and the zero ghost values pull the flow toward zero next to the last row and
	column.
	Table~\ref{tab:ex_tv_3x3} summarizes the key numerical values.

	\begin{table}[h]
		\centering
		\caption{Key numerical values for Example~\ref{ex:tv_3x3}.}
		\label{tab:ex_tv_3x3}
		\begin{tabular}{lc}
			\toprule
			Quantity                                          & Value     \\
			\midrule
			$Q(0)=\|w_0-w_g\|_2^2$                                        & $6.00$    \\
			$\|w_g\|_2^2$                                                  & $8.00$    \\
			$Q(0)<\|w_g\|_2^2$ (improvement over zero flow, \eqref{cond:CLG-TV-Exist})     & Yes       \\
			Condition~\eqref{cond:CLG-TV-NZ1a} satisfied                  & Yes       \\
			$\mathcal{R}_{\mathrm{var}}=0$, $\mathcal{R}_{\mathrm{mix}}=0$ ($M_p=I$) & ---  \\
			Alignment $S_u+S_v$ in~\eqref{cond:CLG-TV-NZ1b}              & $17.91$   \\
			$D^+Q(0)$ (forward difference; $-2(S_u+S_v)\approx-35.81$)             & $-35.75$  \\
			$\alpha^*$                                                     & $0.211$   \\
			$Q(\alpha^*)$                                                  & $2.558$   \\
			Error reduction                                                & $57.4\%$ \\
			Computed $w_\alpha=0$ for all sampled $\alpha\ge$ & $0.72$ \\
			\bottomrule
		\end{tabular}
	\end{table}

	\begin{figure}[!t]
		\centering
		\includegraphics[width=0.62\linewidth]{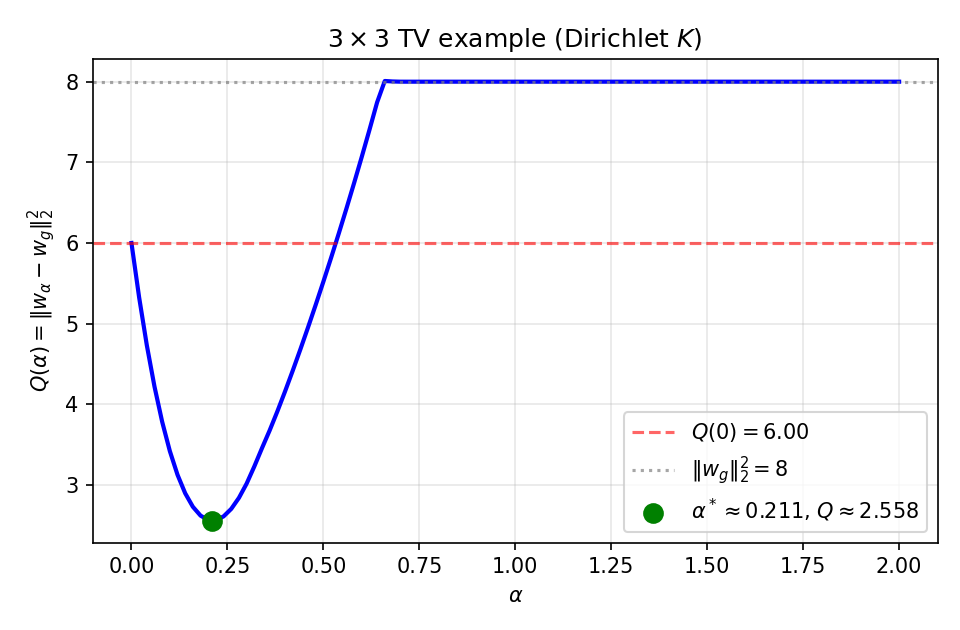}
		\caption{Objective $Q(\alpha) = \|w_\alpha - w_g\|_2^2$ versus $\alpha$ for the $3 \times 3$ TV example. The minimum is at $\alpha^* \approx 0.211$. For $\alpha\ge0.72$ the computed flow is zero and $Q(\alpha)=\|w_g\|_2^2=8$}
		\label{fig:tv_qalpha_ex}
	\end{figure}
\end{example}

\begin{example}
	\label{ex:tv_50x50}
	We consider a $50\times50$ scene containing a centered $10\times10$ square of intensity $1.0$
	on a zero background. The translation is $(0.7,-0.3)$, so that the square moves while the
	background stays static, and the noise has $\sigma=0.10$ and seed $42$. The warp uses zero values
	outside the image, and the noise is generated with NumPy's \texttt{RandomState}. The
	ground-truth flow $w_g$ equals $(0.7,-0.3)$ on the pixels of the square in $I_1$ and zero
	elsewhere, so that $\|w_g\|_2^2=58.0$.
	The accelerated step sizes use the strong convexity modulus
	$\lambda_{\min}(\mathcal{M}_{\mathrm{LK}})\approx1.1\times10^{-2}$ of the data term. The sweep tolerance is $10^{-9}$. The largest
	final relative gap over all sweep solves was $6\times10^{-8}$. The derivative estimates use
	the tolerance $10^{-13}$. Both the derivative estimates and the confirmation run at $\alpha^*$ use at
	most $4\times10^5$ iterations. The parameter is sampled at $51$
	log-spaced points over $[10^{-4},10]$, followed by a $31$-point linear
	refinement in the bracket surrounding the coarse minimum.

	We first check the conditions of Corollary~\ref{thm:tv_exist}. Numerically,
	$Q(0)=\|w_0-w_g\|_2^2 \approx 381.9$.
	Since $Q(0)>\|w_g\|_2^2$, the sufficient condition~\eqref{cond:CLG-TV-Exist}
	\emph{fails}: the unregularized Lucas--Kanade solution $w_0$ fits noise in background
	pixels and produces spurious flow that is farther from $w_g$ than the zero-flow baseline.
	By Corollary~\ref{thm:tv_exist}, an optimal parameter exists whether or not
	condition~\eqref{cond:CLG-TV-Exist} holds. The condition is sufficient, but not necessary, for the
	optimum to improve on the zero flow.

	The computed $w_0$ satisfies $\|\nabla u_0(p)\|_2>0$ and $\|\nabla v_0(p)\|_2>0$ at all $2500$ grid
	points. Hence $\Omega_{11}=\Omega$, and condition~\eqref{cond:CLG-TV-NZ1a} is satisfied.

	Condition~\eqref{cond:CLG-TV-NZ1b} is verified directly with the stabilized
	per-pixel inverse $M_p=(A_p^\top A_p+10^{-2}I)^{-1}$.
	Since $\Omega_{10}=\Omega_{01}=\varnothing$, the mixed-pixel coupling
	$\mathcal R_{\mathrm{mix}}=0$ structurally, independent of the off-diagonal
	entries $m_{12}(p)$.
	The remaining sums evaluate to
	$S_u\approx 24177$, $S_v\approx 24832$, $\mathcal R_{\mathrm{var}}\approx 31638$,
	so that $S_u+S_v-\sqrt2\,\mathcal R_{\mathrm{var}}-\mathcal R_{\mathrm{mix}}\approx 4267>0$.
	Condition~\eqref{cond:CLG-TV-NZ1b} therefore holds, and
	Theorem~\ref{thm:tv_nonzero} guarantees $\alpha^*>0$.
	Forward-difference estimates give $D^+Q(0)\approx-9.6\times10^4$ with step
	$10^{-5}$ and $D^+Q(0)\approx-9.3\times10^4$ with step $10^{-4}$. Both lie below the bound
	\[-2\bigl[(S_u+S_v)-\sqrt2\,\mathcal R_{\mathrm{var}}-\mathcal R_{\mathrm{mix}}\bigr]\approx
	-8.5\times10^3\]
	obtained in the proof of Theorem~\ref{thm:tv_nonzero} (which is consistent with that
	theorem).

	The sweep then locates the minimum at $\alpha^*\approx1.68$ with
	$Q(\alpha^*)\approx29.3$, an error reduction of approximately $92.3\%$.
	A confirmation solve at $\alpha^*$
	reproduces $Q(\alpha^*)$ to four decimal places. For every sampled $\alpha\ge3.2$ the
	computed flow is zero to solver accuracy ($\max_p|w_\alpha(p)|\le5.3\times10^{-12}$), and
	$Q(\alpha)=\|w_g\|_2^2$ to six digits. A dual certificate shows that the exact lower-level
	solution vanishes as well. By~\eqref{LL_OC_TV}, $w_\alpha=0$ if and only if
	$f_{\mathrm{LK}}=-\alpha\,\mathbf{div}\,z$ for a dual field $z$ with $\|z_u(p)\|_2\le1$ and
	$\|z_v(p)\|_2\le1$ at every pixel. A least-norm correction of the solver's dual field at the
	sample $\alpha\approx3.16$ gives a field $y$ with $-\mathbf{div}\,y=f_{\mathrm{LK}}$ (relative residual
	$1.3\times10^{-15}$) and 
	\[\alpha_{\mathrm{cert}}:=\max_p\max\{\|y_u(p)\|_2,\|y_v(p)\|_2\}\approx3.16054.\]
	The field $z=y/\alpha$
	is therefore admissible for every $\alpha\ge\alpha_{\mathrm{cert}}$, and $w_\alpha=0$ for all these $\alpha$. No
	parameter beyond $\alpha_{\mathrm{cert}}$ can therefore improve on the zero flow. The reported $\alpha^*$ gives the
	smallest loss found by the coarse and refined sweeps, and these samples do not certify global
	optimality on $[0,\alpha_{\mathrm{cert}}]$. In
	particular $Q(\alpha^*)<\|w_g\|_2^2$, so the learned regularization improves on the zero
	flow as well as on the unregularized estimate.
	
	TV produces a single block of nearly uniform flow around the moving square. Away from the
	square and a two-pixel band around it, the squared error falls from $259.9$ to $1.7$, and
	$27.6$ of the remaining $29.3$ lies on the square and that band. The block extends a few
	pixels beyond the square, and the motion inside it is underestimated: the horizontal
	component is about $0.40$ instead of $0.7$, and the vertical component is suppressed. This
	loss of contrast is characteristic of TV regularization~\cite{strong2003edge}. On the same scene, quadratic
	regularization with its parameter selected in the same way attains only
	$Q(\alpha^*)\approx43.7$, at $\alpha^*\approx170$. Its mean flow magnitude on the square is
	about $0.19$, and the flow decays only slowly into the background.
	Table~\ref{tab:ex_tv_50x50} summarizes the key numerical values.
	The full curve $Q(\alpha)$ is shown in Figure~\ref{fig:tv_50x50}(d); the pointwise
	error maps before and after TV regularization are displayed in
	Figure~\ref{fig:tv_50x50}(e)--(f).

	\begin{table}[h]
		\centering
		\caption{Key numerical values for Example~\ref{ex:tv_50x50}.}
		\label{tab:ex_tv_50x50}
		\begin{tabular}{lc}
			\toprule
			Quantity & Value \\
			\midrule
			$Q(0)=\|w_0-w_g\|_2^2$ & $381.9$ \\
			$\|w_g\|_2^2$ & $58.0$ \\
			$Q(0)<\|w_g\|_2^2$ (improvement over zero flow, \eqref{cond:CLG-TV-Exist}) & No \\
			Condition~\eqref{cond:CLG-TV-NZ1a} satisfied & Yes \\
			$|\Omega_{11}|/|\Omega_{10}|/|\Omega_{01}|/|\Omega_{00}|$ & $2500/0/0/0$ \\
			Condition~\eqref{cond:CLG-TV-NZ1b} verified & Yes \\
			$S_u$ & $24177$ \\
			$S_v$ & $24832$ \\
			$\mathcal R_{\mathrm{var}}$ & $31638$ \\
			$\mathcal R_{\mathrm{mix}}$ (structural zero, $\Omega_{10}=\Omega_{01}=\varnothing$) & $0$ \\
			$S_u+S_v-\sqrt2\,\mathcal R_{\mathrm{var}}-\mathcal R_{\mathrm{mix}}$ & $4267$ \\
			$D^+Q(0)$ (numerical, step $10^{-5}$) & $-9.6\times 10^4$ \\
			$\alpha^*$ & $1.68$ \\
			$Q(\alpha^*)$ & $29.3$ \\
			Reduction vs.\ $\alpha=0$ & $92.3\%$ \\
			$Q(\alpha^*)<\|w_g\|_2^2$ (improves on zero flow) & Yes \\
			$w_\alpha=0$ (dual certificate) for all $\alpha\ge$ & $3.161$ \\
			\bottomrule
		\end{tabular}
	\end{table}

	\begin{figure}[!htbp]
		\centering
		\includegraphics[width=0.95\textwidth]{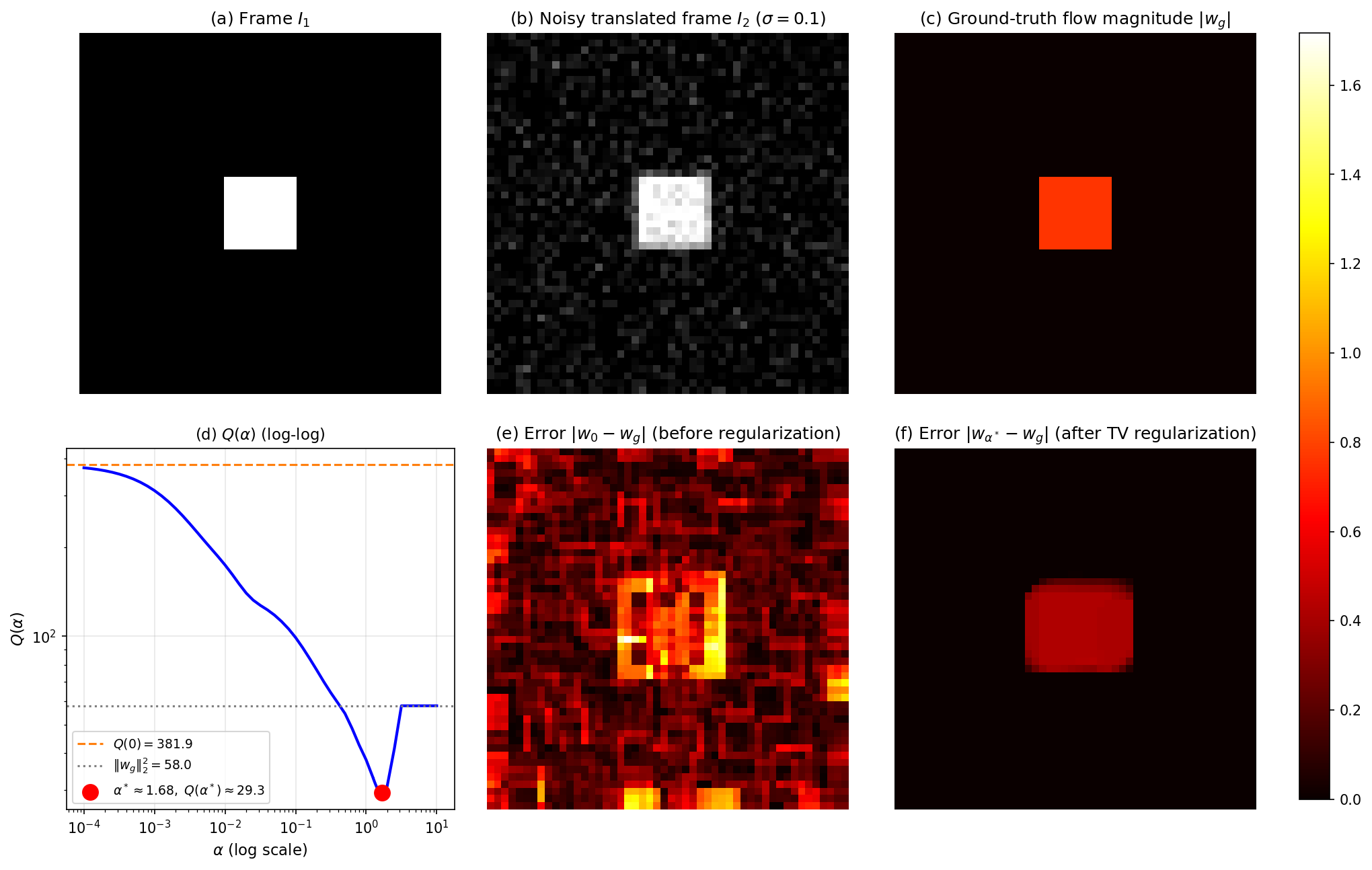}
		\caption{\textbf{Example~\ref{ex:tv_50x50}:}
			(a) Frame $I_1$: centered $10\times10$ square of intensity $1$ (white) on a zero
			(black) background.
			(b) Noisy translated frame $I_2$ ($\sigma=0.10$).
			(c) Ground-truth flow magnitude $|w_g|$: displacement $(0.7,-0.3)$ on the square,
			zero elsewhere (magnitude $\approx0.76$).
			(d) $Q(\alpha)$ on log--log axes: $Q(0)\approx381.9>\|w_g\|_2^2=58.0$ (condition~\eqref{cond:CLG-TV-Exist}
			fails). The minimum $Q(\alpha^*)\approx29.3$ at
			$\alpha^*\approx1.68$ (red marker) lies below $\|w_g\|_2^2$. For $\alpha\ge3.161$ the
			lower-level solution is zero, by a dual certificate, and $Q(\alpha)=\|w_g\|_2^2$.
			(e) Pointwise error $|w_0-w_g|$ before regularization: large spurious
			errors in the static background.
			(f) Pointwise error $|w_{\alpha^*}-w_g|$ after TV regularization ($92.3\%$ reduction
			in $Q$): the background error is largely removed, and the remaining error is a nearly uniform
			block over the square and a narrow band around it, which reflects the loss of contrast
			of TV}
		\label{fig:tv_50x50}
	\end{figure}
\end{example}

The next example is designed to exercise the full condition~\eqref{cond:CLG-TV-NZ1b},
including the mixed-pixel coupling term $\mathcal R_{\mathrm{mix}}$ involving the
off-diagonal entry $m_{12}$.  Examples~\ref{ex:tv_3x3} and~\ref{ex:tv_50x50}
both have $\Omega_{10}=\Omega_{01}=\varnothing$, so $\mathcal R_{\mathrm{mix}}=0$
structurally regardless of $m_{12}$.  Here we construct a $4\times4$ scene in
which $|\Omega_{10}|=9$ and $m_{12}\not\equiv0$, so that $\mathcal R_{\mathrm{mix}}>0$.

\begin{example}
	\label{ex:tv_designed_mix}
	
	We set the inverse structure tensor to be spatially uniform:
	\[
		M_p = (A_p^\top A_p)^{-1} =
		\begin{pmatrix} 1 & -\tfrac{1}{2} \\[2pt] -\tfrac{1}{2} & 1 \end{pmatrix}
		\quad\text{for every pixel }p,
	\]
	so $m_{12}=-\tfrac{1}{2}\neq0$.
	Because $M_p$ is spatially uniform, $\rho(p)=0$ for every pixel and
	$\mathcal R_{\mathrm{var}}=0$ exactly.

	Pixel coordinates are $(i,j)$ with $i,j\in\{0,1,2,3\}$.
	We prescribe piecewise-constant flows with a horizontal step at column $j=2$:
	\[
		u_0(i,j) = \begin{cases} 0 & j \le 1, \\ \tfrac{3}{2} & j \ge 2, \end{cases}
		\quad v_0 \equiv 0;\qquad
		u_g(i,j) = \begin{cases} 0 & j \le 1, \\ 1 & j \ge 2, \end{cases}
		\quad v_g \equiv 0.
	\]
	The residual is $r_u(i,j)=\tfrac{1}{2}$ for $j\ge2$ and zero otherwise; $r_v\equiv0$.
	Hence $Q(0)=\|w_0-w_g\|_2^2=2$ and
	$\|w_g\|_2^2=8$, so $Q(0)<\|w_g\|_2^2$ and
	condition~\eqref{cond:CLG-TV-Exist} is satisfied (Corollary~\ref{thm:tv_exist}).

With the gradient $K$, $\nabla u_0(p)\neq0$ at nine
	pixels. These are the four pixels $(i,1)$, where $\nabla u_0=(\tfrac32,0)^\top$; the three pixels
	$(i,3)$ with $i\le2$, where $\nabla u_0=(-\tfrac32,0)^\top$; the pixel $(3,2)$, where
	$\nabla u_0=(0,-\tfrac32)^\top$; and the pixel $(3,3)$, where $\nabla u_0=(-\tfrac32,-\tfrac32)^\top$. Since
	$v_0\equiv0$, $\nabla v_0\equiv0$ throughout. The partition is $|\Omega_{11}|=0$, $|\Omega_{10}|=9$,
	$|\Omega_{01}|=0$, $|\Omega_{00}|=7$.

	The step in $u_g$ has the same support as $u_0$, so wherever $\nabla u_0(p)=0$
	one also has $\nabla u_g(p)=0$.  Since $v_g\equiv0$, the $v$-implication holds
	trivially.  Condition~\eqref{cond:CLG-TV-NZ1a} is satisfied.

Since $r_u=u_0-u_g=\tfrac13u_0$ and $r_v\equiv0$, we have $\nabla r_u=\tfrac13\nabla u_0$ and
	$\nabla r_v\equiv0$. At every $p\in\Omega_{10}$ the gradient $\nabla r_u(p)$ is therefore aligned with
	$\hat z_u(p)=\nabla u_0(p)/\|\nabla u_0(p)\|_2$, and $\langle\nabla r_u(p),\hat z_u(p)\rangle=\|\nabla r_u(p)\|_2$. This norm equals
	$\tfrac12$ at eight of the nine pixels and $\tfrac{\sqrt2}{2}$ at $(3,3)$. Hence
	\begin{align*}
		S_u &= \sum_{p\in\Omega_{10}}
		\bigl[m_{11}\langle\nabla r_u,\hat z_u\rangle
		+ m_{12}\langle\nabla r_v,\hat z_u\rangle\bigr]
		= 8\cdot\tfrac12+\tfrac{\sqrt2}{2} \approx 4.707, \\
		S_v &= 0 \qquad (\Omega_{11}=\Omega_{01}=\varnothing), \\
		\mathcal R_{\mathrm{var}} &= 0 \qquad (\rho\equiv0), \\
		\mathcal R_{\mathrm{mix}} &= \sum_{p\in\Omega_{10}}|m_{12}(p)|\,\|\nabla r_u(p)\|_2
		= \tfrac12\,S_u \approx 2.354.
	\end{align*}
	Therefore $S_u+S_v-\sqrt2\,\mathcal R_{\mathrm{var}}-\mathcal R_{\mathrm{mix}}=\tfrac12\,S_u\approx2.354>0$, and
	Theorem~\ref{thm:tv_nonzero} guarantees $\alpha^*>0$.
	Table~\ref{tab:ex_tv_designed_mix} summarizes the computed quantities.

	\begin{table}[h]
		\centering
		\caption{Key numerical values for Example~\ref{ex:tv_designed_mix}.}
		\label{tab:ex_tv_designed_mix}
		\begin{tabular}{lc}
			\toprule
			Quantity & Value \\
			\midrule
			$M_p$ (spatially uniform) & $\bigl[\begin{smallmatrix}1&-\frac12\\-\frac12&1\end{smallmatrix}\bigr]$ \\
			$Q(0)=\|w_0-w_g\|_2^2$ & $2$ \\
			$\|w_g\|_2^2$ & $8$ \\
			$Q(0)<\|w_g\|_2^2$ (improvement over zero flow, \eqref{cond:CLG-TV-Exist}) & Yes \\
			Condition~\eqref{cond:CLG-TV-NZ1a} satisfied & Yes \\
			$|\Omega_{11}|/|\Omega_{10}|/|\Omega_{01}|/|\Omega_{00}|$ & $0/9/0/7$ \\
			$S_u$ & $4.707$ \\
			$S_v$ & $0$ \\
			$\mathcal R_{\mathrm{var}}$ & $0$ \\
			$\mathcal R_{\mathrm{mix}}$ ($m_{12}\neq0$, $|\Omega_{10}|=9$) & $2.354$ \\
			$S_u+S_v-\sqrt2\,\mathcal R_{\mathrm{var}}-\mathcal R_{\mathrm{mix}}$ & $2.354$ \\
			Condition~\eqref{cond:CLG-TV-NZ1b} verified & Yes \\
			$\alpha^*>0$ (Theorem~\ref{thm:tv_nonzero}) & Guaranteed \\
			\bottomrule
		\end{tabular}
	\end{table}
\end{example}

\section{Limitations and Conclusion}\label{sec:conclusion}

\subsection{Limitations}\label{subsec:limitations}

The results of this paper hold under assumptions that restrict their scope, and some of the conditions
are conservative. We summarize these restrictions below.

The analysis requires assumption (A1), that is, full-rank local data matrices. This assumption fails in
textureless regions. The $50\times50$ examples therefore use the stabilization $A_p^\top A_p+10^{-2}I$,
which makes $\mathcal{M}_{\mathrm{LK}}$ positive definite, and their results apply to the stabilized data term.

A single scalar parameter is learned from one pair of frames with an exact ground-truth flow. We do not
analyze the effect of errors in the reference flow, nor whether a parameter learned on one pair remains
suitable for other pairs. The examples are small synthetic scenes chosen to illustrate the conditions,
and they do not evaluate performance on benchmark data.

The existence, positivity and second-order conditions are sufficient but not necessary, and they can be
conservative. Condition~\eqref{cond:Exist-Dir} fails for the noisier realizations of
Example~\ref{ex:quad_50x50}, at the tested noise levels $\sigma\ge0.20$ (Remark~\ref{rem:high_noise_quad}
and Figure~\ref{fig:noise_sensitivity}). Condition~\eqref{cond:CLG-TV-Exist}, which guarantees an improvement over the
zero flow, fails in Example~\ref{ex:tv_50x50}. In all these cases an optimal parameter still exists and
improves on $w_0$. Likewise, the residual bound~\eqref{eq:residual_bound} of Theorem~\ref{th:sosc}(b) fails
in Examples~\ref{ex:quad_3x3} and~\ref{ex:quad_50x50}, although $Q''(\alpha^*)>0$ there, and only the
constructed Example~\ref{ex:sosc_constructed} satisfies it.

The results for TV are weaker than those for the quadratic model. We do not derive a second-order
condition for the TV model, and positivity of the optimal parameter is obtained from a Dini-derivative test at $\alpha=0$,
on grids with at least two rows and two columns. Condition~\eqref{cond:CLG-TV-NZ1a} is an idealized
assumption, and $\sqrt2\,\mathcal{R}_{\mathrm{var}}$ and $\mathcal{R}_{\mathrm{mix}}$ are worst-case bounds.

Finally, for neither model do we show that $Q$ has a unique global minimizer. The examples locate the
minimizers numerically, by Brent's method for the quadratic model and by grid refinement for TV, and
these searches do not certify global optimality.

\subsection{Conclusion}

This paper analyzed a bilevel formulation for learning the scalar regularization parameter of CLG
optical flow models from ground-truth flow, in a finite-dimensional setting in which the lower-level
solution is unique.

For quadratic regularization, Proposition~\ref{prop:exist} gives existence under
condition~\eqref{cond:Exist-Dir}, which asks that the unregularized Lucas--Kanade estimate be closer to the
ground truth than the zero flow. Theorem~\ref{th:clg-lossE2} shows that every minimizer, if one exists,
satisfies $\alpha^*>0$ when the first-order
smoothing correction $-\mathcal{M}_{\mathrm{LK}}^{-1}Lw_0$ has a negative component along the residual
$w_0-w_g$. Theorem~\ref{th:sosc} gives second-order sufficient conditions for strict local optimality
through a residual bound derived from a Cauchy--Schwarz estimate. The examples show that this bound is
conservative.

For TV regularization, Corollary~\ref{thm:tv_exist} shows that an optimal parameter always exists,
because the TV solution vanishes beyond a finite value of $\alpha$. Condition~\eqref{cond:CLG-TV-Exist} then
guarantees that the optimum improves on the zero flow. Theorem~\ref{thm:tv_nonzero} ensures $\alpha^*>0$
under two conditions. The structural condition~\eqref{cond:CLG-TV-NZ1a} asks that the ground truth be flat
wherever the unregularized estimate is flat. The quantitative condition~\eqref{cond:CLG-TV-NZ1b} asks that
the weighted alignment $S_u+S_v$ exceed the worst-case bounds $\sqrt2\,\mathcal{R}_{\mathrm{var}}$, for the
spatial variation of $(A_p^\top A_p)^{-1}$, and $\mathcal{R}_{\mathrm{mix}}$, for the pixels at which only one
component of $\nabla w_0$ vanishes. Corollary~\ref{cor:tv_simple} simplifies the condition when the local
structure tensors are spatially constant. In Example~\ref{ex:tv_50x50}, TV regularization
removes most of the background error, and the motion on the square is underestimated.

Several questions remain open. Spatially varying parameters $\alpha(p)$ can adapt the regularization to
the image content. The scalar analysis suggests that such an extension will involve the same
quantities: the residual of the unregularized estimate, its gradients, and the local structure tensors.
Second-order conditions for TV would require a description of the nonsmooth dependence of $w_\alpha$ on
$\alpha$. Smoothing the TV term, for instance with the Huber function~\cite{kalmoun2018investigation}, is one
possible route. The Huber function is only once continuously differentiable, however, so the
second-order argument of Section~\ref{sec:quad} does not carry over directly.

\backmatter

\bmhead{Supplementary information}
The Python code that reproduces every example, table and figure of Section~\ref{sec:examples} is
provided as ancillary files of this arXiv submission (directory \texttt{anc/}). Its README lists, for each example, the script, its inputs, its outputs and
its location in the paper.

\section*{Declarations}

\bmhead{Funding}
This work was supported by the SSE Seed Funding program of Al Akhawayn University in Ifrane (AUI).

\bmhead{Competing interests}
The author has no relevant financial or non-financial interests to disclose.

\bmhead{Data availability}
All data are synthetic and are generated by the code in the ancillary files. No external data
sets were used.

\bmhead{Code availability}
The code is provided as ancillary files of this arXiv submission.

\bibliography{references}


\begin{thebibliography}{44}
\ifx \bisbn   \undefined \def \bisbn  #1{ISBN #1}\fi
\ifx \binits  \undefined \def \binits#1{#1}\fi
\ifx \bauthor  \undefined \def \bauthor#1{#1}\fi
\ifx \batitle  \undefined \def \batitle#1{#1}\fi
\ifx \bjtitle  \undefined \def \bjtitle#1{#1}\fi
\ifx \bvolume  \undefined \def \bvolume#1{\textbf{#1}}\fi
\ifx \byear  \undefined \def \byear#1{#1}\fi
\ifx \bissue  \undefined \def \bissue#1{#1}\fi
\ifx \bfpage  \undefined \def \bfpage#1{#1}\fi
\ifx \blpage  \undefined \def \blpage #1{#1}\fi
\ifx \burl  \undefined \def \burl#1{\textsf{#1}}\fi
\ifx \doiurl  \undefined \def \doiurl#1{\url{https://doi.org/#1}}\fi
\ifx \betal  \undefined \def \betal{\textit{et al.}}\fi
\ifx \binstitute  \undefined \def \binstitute#1{#1}\fi
\ifx \binstitutionaled  \undefined \def \binstitutionaled#1{#1}\fi
\ifx \bctitle  \undefined \def \bctitle#1{#1}\fi
\ifx \beditor  \undefined \def \beditor#1{#1}\fi
\ifx \bpublisher  \undefined \def \bpublisher#1{#1}\fi
\ifx \bbtitle  \undefined \def \bbtitle#1{#1}\fi
\ifx \bedition  \undefined \def \bedition#1{#1}\fi
\ifx \bseriesno  \undefined \def \bseriesno#1{#1}\fi
\ifx \blocation  \undefined \def \blocation#1{#1}\fi
\ifx \bsertitle  \undefined \def \bsertitle#1{#1}\fi
\ifx \bsnm \undefined \def \bsnm#1{#1}\fi
\ifx \bsuffix \undefined \def \bsuffix#1{#1}\fi
\ifx \bparticle \undefined \def \bparticle#1{#1}\fi
\ifx \barticle \undefined \def \barticle#1{#1}\fi
\bibcommenthead
\ifx \bconfdate \undefined \def \bconfdate #1{#1}\fi
\ifx \botherref \undefined \def \botherref #1{#1}\fi
\ifx \url \undefined \def \url#1{\textsf{#1}}\fi
\ifx \bchapter \undefined \def \bchapter#1{#1}\fi
\ifx \bbook \undefined \def \bbook#1{#1}\fi
\ifx \bcomment \undefined \def \bcomment#1{#1}\fi
\ifx \oauthor \undefined \def \oauthor#1{#1}\fi
\ifx \citeauthoryear \undefined \def \citeauthoryear#1{#1}\fi
\ifx \endbibitem  \undefined \def \endbibitem {}\fi
\ifx \bconflocation  \undefined \def \bconflocation#1{#1}\fi
\ifx \arxivurl  \undefined \def \arxivurl#1{\textsf{#1}}\fi
\csname PreBibitemsHook\endcsname

\bibitem[\protect\citeauthoryear{Lucas and Kanade}{1981}]{lucas1981iterative}
\begin{bchapter}
\bauthor{\bsnm{Lucas}, \binits{B.D.}},
\bauthor{\bsnm{Kanade}, \binits{T.}}:
\bctitle{An iterative image registration technique with an application to
  stereo vision}.
In: \bbtitle{Proceedings of the 7th International Joint Conference on
  Artificial Intelligence (IJCAI)},
vol. \bseriesno{2},
pp. \bfpage{674}--\blpage{679}
(\byear{1981})
\end{bchapter}
\endbibitem

\bibitem[\protect\citeauthoryear{Horn and Schunck}{1981}]{horn1981determining}
\begin{barticle}
\bauthor{\bsnm{Horn}, \binits{B.K.P.}},
\bauthor{\bsnm{Schunck}, \binits{B.G.}}:
\batitle{Determining optical flow}.
\bjtitle{Artif. Intell.}
\bvolume{17}(\bissue{1--3}),
\bfpage{185}--\blpage{203}
(\byear{1981})
\doiurl{10.1016/0004-3702(81)90024-2}
\end{barticle}
\endbibitem

\bibitem[\protect\citeauthoryear{Rudin et~al.}{1992}]{rudin1992nonlinear}
\begin{barticle}
\bauthor{\bsnm{Rudin}, \binits{L.I.}},
\bauthor{\bsnm{Osher}, \binits{S.}},
\bauthor{\bsnm{Fatemi}, \binits{E.}}:
\batitle{Nonlinear total variation based noise removal algorithms}.
\bjtitle{Physica D}
\bvolume{60}(\bissue{1--4}),
\bfpage{259}--\blpage{268}
(\byear{1992})
\doiurl{10.1016/0167-2789(92)90242-F}
\end{barticle}
\endbibitem

\bibitem[\protect\citeauthoryear{Shulman and
  Herv{\'e}}{1989}]{shulman1989regularization}
\begin{bchapter}
\bauthor{\bsnm{Shulman}, \binits{D.}},
\bauthor{\bsnm{Herv{\'e}}, \binits{J.-Y.}}:
\bctitle{Regularization of discontinuous flow fields}.
In: \bbtitle{Proceedings of the IEEE Workshop on Visual Motion},
\bconflocation{Irvine, CA, USA},
pp. \bfpage{81}--\blpage{86}
(\byear{1989}).
\doiurl{10.1109/WVM.1989.47097}
\end{bchapter}
\endbibitem

\bibitem[\protect\citeauthoryear{Rouchouze et~al.}{1994}]{rouchouze1994motion}
\begin{bchapter}
\bauthor{\bsnm{Rouchouze}, \binits{B.}},
\bauthor{\bsnm{Mathieu}, \binits{P.}},
\bauthor{\bsnm{Gaidon}, \binits{T.}},
\bauthor{\bsnm{Barlaud}, \binits{M.}}:
\bctitle{Motion estimation based on {Markov} random fields}.
In: \bbtitle{Proceedings of the 1st IEEE International Conference on Image
  Processing, Austin, TX, USA},
vol. \bseriesno{3},
pp. \bfpage{270}--\blpage{274}
(\byear{1994}).
\doiurl{10.1109/ICIP.1994.413845}
\end{bchapter}
\endbibitem

\bibitem[\protect\citeauthoryear{Deriche et~al.}{1996}]{deriche1995optical}
\begin{bchapter}
\bauthor{\bsnm{Deriche}, \binits{R.}},
\bauthor{\bsnm{Kornprobst}, \binits{P.}},
\bauthor{\bsnm{Aubert}, \binits{G.}}:
\bctitle{Optical-flow estimation while preserving its discontinuities: A
  variational approach}.
In: \beditor{\bsnm{Li}, \binits{S.Z.}},
\beditor{\bsnm{Mital}, \binits{D.P.}},
\beditor{\bsnm{Teoh}, \binits{E.K.}},
\beditor{\bsnm{Wang}, \binits{H.}} (eds.)
\bbtitle{Recent Developments in Computer Vision (ACCV'95)}.
\bsertitle{Lecture Notes in Computer Science},
vol. \bseriesno{1035},
pp. \bfpage{69}--\blpage{80}.
\bpublisher{Springer},
\blocation{Berlin, Heidelberg}
(\byear{1996}).
\doiurl{10.1007/3-540-60793-5_63}
\end{bchapter}
\endbibitem

\bibitem[\protect\citeauthoryear{Zach et~al.}{2007}]{zach2007aduality}
\begin{bchapter}
\bauthor{\bsnm{Zach}, \binits{C.}},
\bauthor{\bsnm{Pock}, \binits{T.}},
\bauthor{\bsnm{Bischof}, \binits{H.}}:
\bctitle{A duality based approach for realtime {TV-L1} optical flow}.
In: \beditor{\bsnm{Hamprecht}, \binits{F.A.}},
\beditor{\bsnm{Schn{\"o}rr}, \binits{C.}},
\beditor{\bsnm{J{\"a}hne}, \binits{B.}} (eds.)
\bbtitle{Pattern Recognition (DAGM 2007)}.
\bsertitle{Lecture Notes in Computer Science},
vol. \bseriesno{4713},
pp. \bfpage{214}--\blpage{223}.
\bpublisher{Springer},
\blocation{Berlin, Heidelberg}
(\byear{2007}).
\doiurl{10.1007/978-3-540-74936-3_22}
\end{bchapter}
\endbibitem

\bibitem[\protect\citeauthoryear{Wedel et~al.}{2009}]{Wedel2008}
\begin{bchapter}
\bauthor{\bsnm{Wedel}, \binits{A.}},
\bauthor{\bsnm{Pock}, \binits{T.}},
\bauthor{\bsnm{Zach}, \binits{C.}},
\bauthor{\bsnm{Bischof}, \binits{H.}},
\bauthor{\bsnm{Cremers}, \binits{D.}}:
\bctitle{An improved algorithm for {TV-L1} optical flow}.
In: \beditor{\bsnm{Cremers}, \binits{D.}},
\beditor{\bsnm{Rosenhahn}, \binits{B.}},
\beditor{\bsnm{Yuille}, \binits{A.L.}},
\beditor{\bsnm{Schmidt}, \binits{F.R.}} (eds.)
\bbtitle{Statistical and Geometrical Approaches to Visual Motion Analysis}.
\bsertitle{Lecture Notes in Computer Science},
vol. \bseriesno{5604},
pp. \bfpage{23}--\blpage{45}.
\bpublisher{Springer},
\blocation{Berlin, Heidelberg}
(\byear{2009}).
\doiurl{10.1007/978-3-642-03061-1_2}
\end{bchapter}
\endbibitem

\bibitem[\protect\citeauthoryear{Chambolle and
  Pock}{2011}]{chambolle2011afirst}
\begin{barticle}
\bauthor{\bsnm{Chambolle}, \binits{A.}},
\bauthor{\bsnm{Pock}, \binits{T.}}:
\batitle{A first-order primal-dual algorithm for convex problems with
  applications to imaging}.
\bjtitle{J. Math. Imaging Vis.}
\bvolume{40}(\bissue{1}),
\bfpage{120}--\blpage{145}
(\byear{2011})
\doiurl{10.1007/s10851-010-0251-1}
\end{barticle}
\endbibitem

\bibitem[\protect\citeauthoryear{Ranftl et~al.}{2014}]{rantfl2014nonlocal}
\begin{bchapter}
\bauthor{\bsnm{Ranftl}, \binits{R.}},
\bauthor{\bsnm{Bredies}, \binits{K.}},
\bauthor{\bsnm{Pock}, \binits{T.}}:
\bctitle{Non-local total generalized variation for optical flow estimation}.
In: \beditor{\bsnm{Fleet}, \binits{D.}},
\beditor{\bsnm{Pajdla}, \binits{T.}},
\beditor{\bsnm{Schiele}, \binits{B.}},
\beditor{\bsnm{Tuytelaars}, \binits{T.}} (eds.)
\bbtitle{Computer Vision -- ECCV 2014}.
\bsertitle{Lecture Notes in Computer Science},
vol. \bseriesno{8689},
pp. \bfpage{439}--\blpage{454}.
\bpublisher{Springer},
\blocation{Cham}
(\byear{2014}).
\doiurl{10.1007/978-3-319-10590-1_29}
\end{bchapter}
\endbibitem

\bibitem[\protect\citeauthoryear{Ali et~al.}{2016}]{ali2016illumination}
\begin{barticle}
\bauthor{\bsnm{Ali}, \binits{S.}},
\bauthor{\bsnm{Daul}, \binits{C.}},
\bauthor{\bsnm{Galbrun}, \binits{E.}},
\bauthor{\bsnm{Blondel}, \binits{W.}}:
\batitle{Illumination invariant optical flow using neighborhood descriptors}.
\bjtitle{Comput. Vis. Image Underst.}
\bvolume{145},
\bfpage{95}--\blpage{110}
(\byear{2016})
\doiurl{10.1016/j.cviu.2015.12.003}
\end{barticle}
\endbibitem

\bibitem[\protect\citeauthoryear{Chambolle}{2004}]{chambolle2004algorithm}
\begin{barticle}
\bauthor{\bsnm{Chambolle}, \binits{A.}}:
\batitle{An algorithm for total variation minimization and applications}.
\bjtitle{J. Math. Imaging Vis.}
\bvolume{20}(\bissue{1--2}),
\bfpage{89}--\blpage{97}
(\byear{2004})
\doiurl{10.1023/B:JMIV.0000011325.36760.1e}
\end{barticle}
\endbibitem

\bibitem[\protect\citeauthoryear{Kalmoun}{2018}]{kalmoun2018investigation}
\begin{barticle}
\bauthor{\bsnm{Kalmoun}, \binits{E.M.}}:
\batitle{An investigation of smooth {TV}-like regularization in the context of
  the optical flow problem}.
\bjtitle{J. Imaging}
\bvolume{4}(\bissue{2}),
\bfpage{31}
(\byear{2018})
\doiurl{10.3390/jimaging4020031}
\end{barticle}
\endbibitem

\bibitem[\protect\citeauthoryear{Bruhn et~al.}{2005}]{bruhn2005lucas}
\begin{barticle}
\bauthor{\bsnm{Bruhn}, \binits{A.}},
\bauthor{\bsnm{Weickert}, \binits{J.}},
\bauthor{\bsnm{Schn{\"o}rr}, \binits{C.}}:
\batitle{{Lucas/Kanade} meets {Horn/Schunck}: Combining local and global optic
  flow methods}.
\bjtitle{Int. J. Comput. Vis.}
\bvolume{61}(\bissue{3}),
\bfpage{211}--\blpage{231}
(\byear{2005})
\doiurl{10.1023/B:VISI.0000045324.43199.43}
\end{barticle}
\endbibitem

\bibitem[\protect\citeauthoryear{Kalmoun and Garrido}{2013}]{kalmoun2013trust}
\begin{barticle}
\bauthor{\bsnm{Kalmoun}, \binits{E.M.}},
\bauthor{\bsnm{Garrido}, \binits{L.}}:
\batitle{Trust region versus line search for computing the optical flow}.
\bjtitle{Multiscale Model. Simul.}
\bvolume{11}(\bissue{3}),
\bfpage{890}--\blpage{906}
(\byear{2013})
\doiurl{10.1137/120897213}
\end{barticle}
\endbibitem

\bibitem[\protect\citeauthoryear{Wedel et~al.}{2009}]{wedel2009structure}
\begin{bchapter}
\bauthor{\bsnm{Wedel}, \binits{A.}},
\bauthor{\bsnm{Cremers}, \binits{D.}},
\bauthor{\bsnm{Pock}, \binits{T.}},
\bauthor{\bsnm{Bischof}, \binits{H.}}:
\bctitle{Structure- and motion-adaptive regularization for high accuracy optic
  flow}.
In: \bbtitle{2009 IEEE 12th International Conference on Computer Vision (ICCV),
  Kyoto, Japan},
pp. \bfpage{1663}--\blpage{1668}
(\byear{2009}).
\doiurl{10.1109/ICCV.2009.5459375}
\end{bchapter}
\endbibitem

\bibitem[\protect\citeauthoryear{Belhachmi and Hecht}{2011}]{Belhachmi2011}
\begin{barticle}
\bauthor{\bsnm{Belhachmi}, \binits{Z.}},
\bauthor{\bsnm{Hecht}, \binits{F.}}:
\batitle{Control of the effects of regularization on variational optic flow
  computations}.
\bjtitle{J. Math. Imaging Vis.}
\bvolume{40}(\bissue{1}),
\bfpage{1}--\blpage{19}
(\byear{2011})
\doiurl{10.1007/s10851-010-0239-x}
\end{barticle}
\endbibitem

\bibitem[\protect\citeauthoryear{Belhachmi and Hecht}{2016}]{Belhachmi2016}
\begin{barticle}
\bauthor{\bsnm{Belhachmi}, \binits{Z.}},
\bauthor{\bsnm{Hecht}, \binits{F.}}:
\batitle{An adaptive approach for the segmentation and the {TV}-filtering in
  the optic flow estimation}.
\bjtitle{J. Math. Imaging Vis.}
\bvolume{54}(\bissue{3}),
\bfpage{358}--\blpage{377}
(\byear{2016})
\doiurl{10.1007/s10851-015-0608-6}
\end{barticle}
\endbibitem

\bibitem[\protect\citeauthoryear{Engl et~al.}{1996}]{engl1996regularization}
\begin{bbook}
\bauthor{\bsnm{Engl}, \binits{H.W.}},
\bauthor{\bsnm{Hanke}, \binits{M.}},
\bauthor{\bsnm{Neubauer}, \binits{A.}}:
\bbtitle{Regularization of Inverse Problems}.
\bsertitle{Mathematics and Its Applications},
vol. \bseriesno{375}.
\bpublisher{Kluwer Academic Publishers},
\blocation{Dordrecht}
(\byear{1996}).
\doiurl{10.1007/978-94-009-1740-8}
\end{bbook}
\endbibitem

\bibitem[\protect\citeauthoryear{Hansen}{1992}]{hansen1992lcurve}
\begin{barticle}
\bauthor{\bsnm{Hansen}, \binits{P.C.}}:
\batitle{Analysis of discrete ill-posed problems by means of the {L}-curve}.
\bjtitle{SIAM Rev.}
\bvolume{34}(\bissue{4}),
\bfpage{561}--\blpage{580}
(\byear{1992})
\doiurl{10.1137/1034115}
\end{barticle}
\endbibitem

\bibitem[\protect\citeauthoryear{Golub et~al.}{1979}]{golub1979gcv}
\begin{barticle}
\bauthor{\bsnm{Golub}, \binits{G.H.}},
\bauthor{\bsnm{Heath}, \binits{M.}},
\bauthor{\bsnm{Wahba}, \binits{G.}}:
\batitle{Generalized cross-validation as a method for choosing a good ridge
  parameter}.
\bjtitle{Technometrics}
\bvolume{21}(\bissue{2}),
\bfpage{215}--\blpage{223}
(\byear{1979})
\doiurl{10.1080/00401706.1979.10489751}
\end{barticle}
\endbibitem

\bibitem[\protect\citeauthoryear{Stein}{1981}]{stein1981sure}
\begin{barticle}
\bauthor{\bsnm{Stein}, \binits{C.M.}}:
\batitle{Estimation of the mean of a multivariate normal distribution}.
\bjtitle{Ann. Stat.}
\bvolume{9}(\bissue{6}),
\bfpage{1135}--\blpage{1151}
(\byear{1981})
\doiurl{10.1214/aos/1176345632}
\end{barticle}
\endbibitem

\bibitem[\protect\citeauthoryear{Baker et~al.}{2011}]{baker2011middlebury}
\begin{barticle}
\bauthor{\bsnm{Baker}, \binits{S.}},
\bauthor{\bsnm{Scharstein}, \binits{D.}},
\bauthor{\bsnm{Lewis}, \binits{J.P.}},
\bauthor{\bsnm{Roth}, \binits{S.}},
\bauthor{\bsnm{Black}, \binits{M.J.}},
\bauthor{\bsnm{Szeliski}, \binits{R.}}:
\batitle{A database and evaluation methodology for optical flow}.
\bjtitle{Int. J. Comput. Vis.}
\bvolume{92}(\bissue{1}),
\bfpage{1}--\blpage{31}
(\byear{2011})
\doiurl{10.1007/s11263-010-0390-2}
\end{barticle}
\endbibitem

\bibitem[\protect\citeauthoryear{Barron et~al.}{1994}]{barron1994performance}
\begin{barticle}
\bauthor{\bsnm{Barron}, \binits{J.L.}},
\bauthor{\bsnm{Fleet}, \binits{D.J.}},
\bauthor{\bsnm{Beauchemin}, \binits{S.S.}}:
\batitle{Performance of optical flow techniques}.
\bjtitle{Int. J. Comput. Vis.}
\bvolume{12}(\bissue{1}),
\bfpage{43}--\blpage{77}
(\byear{1994})
\doiurl{10.1007/BF01420984}
\end{barticle}
\endbibitem

\bibitem[\protect\citeauthoryear{Pagliari et~al.}{2022}]{pagliari2022tvtype}
\begin{barticle}
\bauthor{\bsnm{Pagliari}, \binits{V.}},
\bauthor{\bsnm{Papafitsoros}, \binits{K.}},
\bauthor{\bsnm{Rai{\c{t}}{\u{a}}}, \binits{B.}},
\bauthor{\bsnm{Vikelis}, \binits{A.}}:
\batitle{Bilevel training schemes in imaging for total variation--type
  functionals with convex integrands}.
\bjtitle{SIAM J. Imaging Sci.}
\bvolume{15}(\bissue{4}),
\bfpage{1690}--\blpage{1728}
(\byear{2022})
\doiurl{10.1137/21M1467328}
\end{barticle}
\endbibitem

\bibitem[\protect\citeauthoryear{Kunisch and Pock}{2013}]{kunisch2013bilevel}
\begin{barticle}
\bauthor{\bsnm{Kunisch}, \binits{K.}},
\bauthor{\bsnm{Pock}, \binits{T.}}:
\batitle{A bilevel optimization approach for parameter learning in variational
  models}.
\bjtitle{SIAM J. Imaging Sci.}
\bvolume{6}(\bissue{2}),
\bfpage{938}--\blpage{983}
(\byear{2013})
\doiurl{10.1137/120882706}
\end{barticle}
\endbibitem

\bibitem[\protect\citeauthoryear{{De los Reyes} and
  Sch{\"o}nlieb}{2013}]{delosreyes2013}
\begin{barticle}
\bauthor{\bsnm{{De los Reyes}}, \binits{J.C.}},
\bauthor{\bsnm{Sch{\"o}nlieb}, \binits{C.-B.}}:
\batitle{Image denoising: learning the noise model via nonsmooth
  {PDE}-constrained optimization}.
\bjtitle{Inverse Probl. Imaging}
\bvolume{7}(\bissue{4}),
\bfpage{1183}--\blpage{1214}
(\byear{2013})
\doiurl{10.3934/ipi.2013.7.1183}
\end{barticle}
\endbibitem

\bibitem[\protect\citeauthoryear{{De los Reyes}
  et~al.}{2016}]{delosreyes2016structure}
\begin{barticle}
\bauthor{\bsnm{{De los Reyes}}, \binits{J.C.}},
\bauthor{\bsnm{Sch{\"o}nlieb}, \binits{C.-B.}},
\bauthor{\bsnm{Valkonen}, \binits{T.}}:
\batitle{The structure of optimal parameters for image restoration problems}.
\bjtitle{J. Math. Anal. Appl.}
\bvolume{434}(\bissue{1}),
\bfpage{464}--\blpage{500}
(\byear{2016})
\doiurl{10.1016/j.jmaa.2015.09.023}
\end{barticle}
\endbibitem

\bibitem[\protect\citeauthoryear{{De los Reyes} and
  Villac{\'i}s}{2022}]{delosreyes2022optimality}
\begin{barticle}
\bauthor{\bsnm{{De los Reyes}}, \binits{J.C.}},
\bauthor{\bsnm{Villac{\'i}s}, \binits{D.}}:
\batitle{Optimality conditions for bilevel imaging learning problems with total
  variation regularization}.
\bjtitle{SIAM J. Imaging Sci.}
\bvolume{15}(\bissue{4}),
\bfpage{1646}--\blpage{1689}
(\byear{2022})
\doiurl{10.1137/21M143412X}
\end{barticle}
\endbibitem

\bibitem[\protect\citeauthoryear{{De los Reyes}}{2023}]{delosreyes2023mpcc}
\begin{barticle}
\bauthor{\bsnm{{De los Reyes}}, \binits{J.C.}}:
\batitle{Bilevel imaging learning problems as mathematical programs with
  complementarity constraints: Reformulation and theory}.
\bjtitle{SIAM J. Imaging Sci.}
\bvolume{16}(\bissue{3}),
\bfpage{1655}--\blpage{1686}
(\byear{2023})
\doiurl{10.1137/21M1450744}
\end{barticle}
\endbibitem

\bibitem[\protect\citeauthoryear{Ehrhardt et~al.}{2025}]{ehrhardt2024optimal}
\begin{bchapter}
\bauthor{\bsnm{Ehrhardt}, \binits{M.J.}},
\bauthor{\bsnm{Gazzola}, \binits{S.}},
\bauthor{\bsnm{Scott}, \binits{S.J.}}:
\bctitle{On optimal regularization parameters via bilevel learning}.
In: \beditor{\bsnm{Bubba}, \binits{T.A.}} (ed.)
\bbtitle{Data-driven Models in Inverse Problems}.
\bsertitle{Radon Series on Computational and Applied Mathematics},
vol. \bseriesno{31},
pp. \bfpage{1}--\blpage{38}.
\bpublisher{De Gruyter},
\blocation{Berlin, Boston}
(\byear{2025}).
\doiurl{10.1515/9783111251233-001}
\end{bchapter}
\endbibitem

\bibitem[\protect\citeauthoryear{Hinterm{\"u}ller and
  Wu}{2015}]{hintermuller2015bilevel}
\begin{barticle}
\bauthor{\bsnm{Hinterm{\"u}ller}, \binits{M.}},
\bauthor{\bsnm{Wu}, \binits{T.}}:
\batitle{Bilevel optimization for calibrating point spread functions in blind
  deconvolution}.
\bjtitle{Inverse Probl. Imaging}
\bvolume{9}(\bissue{4}),
\bfpage{1139}--\blpage{1169}
(\byear{2015})
\doiurl{10.3934/ipi.2015.9.1139}
\end{barticle}
\endbibitem

\bibitem[\protect\citeauthoryear{Bartels and
  Weber}{2023}]{bartels2020fractional}
\begin{barticle}
\bauthor{\bsnm{Bartels}, \binits{S.}},
\bauthor{\bsnm{Weber}, \binits{N.}}:
\batitle{Parameter learning and fractional differential operators: Applications
  in regularized image denoising and decomposition problems}.
\bjtitle{Math. Control Relat. Fields}
\bvolume{13}(\bissue{1}),
\bfpage{35}--\blpage{62}
(\byear{2023})
\doiurl{10.3934/mcrf.2021048}
\end{barticle}
\endbibitem

\bibitem[\protect\citeauthoryear{Crockett and
  Fessler}{2022}]{crockett2022bilevel}
\begin{barticle}
\bauthor{\bsnm{Crockett}, \binits{C.}},
\bauthor{\bsnm{Fessler}, \binits{J.A.}}:
\batitle{Bilevel methods for image reconstruction}.
\bjtitle{Found. Trends Signal Process.}
\bvolume{15}(\bissue{2--3}),
\bfpage{121}--\blpage{289}
(\byear{2022})
\doiurl{10.1561/2000000111}
\end{barticle}
\endbibitem

\bibitem[\protect\citeauthoryear{Ehrhardt and
  Roberts}{2021}]{ehrhardt2021dfbilevel}
\begin{barticle}
\bauthor{\bsnm{Ehrhardt}, \binits{M.J.}},
\bauthor{\bsnm{Roberts}, \binits{L.}}:
\batitle{Inexact derivative-free optimization for bilevel learning}.
\bjtitle{J. Math. Imaging Vis.}
\bvolume{63}(\bissue{5}),
\bfpage{580}--\blpage{600}
(\byear{2021})
\doiurl{10.1007/s10851-021-01020-8}
\end{barticle}
\endbibitem

\bibitem[\protect\citeauthoryear{Ehrhardt and
  Roberts}{2024}]{ehrhardt2024hypergradients}
\begin{barticle}
\bauthor{\bsnm{Ehrhardt}, \binits{M.J.}},
\bauthor{\bsnm{Roberts}, \binits{L.}}:
\batitle{Analyzing inexact hypergradients for bilevel learning}.
\bjtitle{IMA J. Appl. Math.}
\bvolume{89}(\bissue{1}),
\bfpage{254}--\blpage{278}
(\byear{2024})
\doiurl{10.1093/imamat/hxad035}
\end{barticle}
\endbibitem

\bibitem[\protect\citeauthoryear{Suonper{\"a} and
  Valkonen}{2025}]{suonpera2025singleloop}
\begin{barticle}
\bauthor{\bsnm{Suonper{\"a}}, \binits{E.}},
\bauthor{\bsnm{Valkonen}, \binits{T.}}:
\batitle{Single-loop methods for bilevel parameter learning in inverse
  imaging}.
\bjtitle{J. Nonsmooth Anal. Optim.}
\bvolume{5},
\bfpage{15577}
(\byear{2025})
\doiurl{10.46298/jnsao-2025-15577}
\end{barticle}
\endbibitem

\bibitem[\protect\citeauthoryear{Domke}{2012}]{domke2012generic}
\begin{bchapter}
\bauthor{\bsnm{Domke}, \binits{J.}}:
\bctitle{Generic methods for optimization-based modeling}.
In: \beditor{\bsnm{Lawrence}, \binits{N.D.}},
\beditor{\bsnm{Girolami}, \binits{M.}} (eds.)
\bbtitle{Proceedings of the Fifteenth International Conference on Artificial
  Intelligence and Statistics (AISTATS)}.
\bsertitle{Proceedings of Machine Learning Research},
vol. \bseriesno{22},
pp. \bfpage{318}--\blpage{326}
(\byear{2012}).
\burl{https://proceedings.mlr.press/v22/domke12.html}
\end{bchapter}
\endbibitem

\bibitem[\protect\citeauthoryear{Domke}{2013}]{domke2013approxmarg}
\begin{barticle}
\bauthor{\bsnm{Domke}, \binits{J.}}:
\batitle{Learning graphical model parameters with approximate marginal
  inference}.
\bjtitle{IEEE Trans. Pattern Anal. Mach. Intell.}
\bvolume{35}(\bissue{10}),
\bfpage{2454}--\blpage{2467}
(\byear{2013})
\doiurl{10.1109/TPAMI.2013.31}
\end{barticle}
\endbibitem

\bibitem[\protect\citeauthoryear{Ranftl and Pock}{2014}]{ranftl2014}
\begin{bchapter}
\bauthor{\bsnm{Ranftl}, \binits{R.}},
\bauthor{\bsnm{Pock}, \binits{T.}}:
\bctitle{A deep variational model for image segmentation}.
In: \beditor{\bsnm{Jiang}, \binits{X.}},
\beditor{\bsnm{Hornegger}, \binits{J.}},
\beditor{\bsnm{Koch}, \binits{R.}} (eds.)
\bbtitle{Pattern Recognition (GCPR 2014)}.
\bsertitle{Lecture Notes in Computer Science},
vol. \bseriesno{8753},
pp. \bfpage{107}--\blpage{118}.
\bpublisher{Springer},
\blocation{Cham}
(\byear{2014}).
\doiurl{10.1007/978-3-319-11752-2_9}
\end{bchapter}
\endbibitem

\bibitem[\protect\citeauthoryear{Jiang et~al.}{2020}]{jiang2020joint}
\begin{bchapter}
\bauthor{\bsnm{Jiang}, \binits{S.}},
\bauthor{\bsnm{Campbell}, \binits{D.}},
\bauthor{\bsnm{Liu}, \binits{M.}},
\bauthor{\bsnm{Gould}, \binits{S.}},
\bauthor{\bsnm{Hartley}, \binits{R.}}:
\bctitle{Joint unsupervised learning of optical flow and egomotion with
  bi-level optimization}.
In: \bbtitle{2020 International Conference on 3D Vision (3DV)},
pp. \bfpage{682}--\blpage{691}
(\byear{2020}).
\doiurl{10.1109/3DV50981.2020.00078}
\end{bchapter}
\endbibitem

\bibitem[\protect\citeauthoryear{El~Hakoume et~al.}{2023}]{elhakoume2023lame}
\begin{barticle}
\bauthor{\bsnm{El~Hakoume}, \binits{A.}},
\bauthor{\bsnm{Laghrib}, \binits{A.}},
\bauthor{\bsnm{Hadri}, \binits{A.}},
\bauthor{\bsnm{Afraites}, \binits{L.}}:
\batitle{An optimal fluid optical flow registration for super-resolution with
  {L}am\'{e} parameters learning}.
\bjtitle{J. Optim. Theory Appl.}
\bvolume{197}(\bissue{2}),
\bfpage{508}--\blpage{538}
(\byear{2023})
\doiurl{10.1007/s10957-023-02186-4}
\end{barticle}
\endbibitem

\bibitem[\protect\citeauthoryear{El~Malki et~al.}{2024}]{elmalki2024nonlocal}
\begin{barticle}
\bauthor{\bsnm{El~Malki}, \binits{I.}},
\bauthor{\bsnm{Jauberteau}, \binits{F.}},
\bauthor{\bsnm{Laghrib}, \binits{A.}},
\bauthor{\bsnm{Nachaoui}, \binits{M.}}:
\batitle{A bilevel learning approach for nonlocal image deblurring with
  variable weights parameter}.
\bjtitle{Digit. Signal Process.}
\bvolume{149},
\bfpage{104505}
(\byear{2024})
\doiurl{10.1016/j.dsp.2024.104505}
\end{barticle}
\endbibitem

\bibitem[\protect\citeauthoryear{Strong and Chan}{2003}]{strong2003edge}
\begin{barticle}
\bauthor{\bsnm{Strong}, \binits{D.}},
\bauthor{\bsnm{Chan}, \binits{T.}}:
\batitle{Edge-preserving and scale-dependent properties of total variation
  regularization}.
\bjtitle{Inverse Probl.}
\bvolume{19}(\bissue{6}),
\bfpage{S165}--\blpage{S187}
(\byear{2003})
\doiurl{10.1088/0266-5611/19/6/059}
\end{barticle}
\endbibitem

\end{thebibliography}

\end{document}